\documentclass[10pt,a4paper]{article}
\usepackage[utf8]{inputenc}
\usepackage[english]{babel}
\usepackage{amsmath}
\usepackage{amsfonts}
\usepackage{amssymb}
\usepackage{graphicx}
\usepackage{tikz,tikz-3dplot}
\usepackage{wrapfig}
\usepackage{caption}
\usepackage{subcaption}
\usepackage{authblk}

\DeclareMathOperator{\const}{const}

\usepackage{xcolor}

\title{A study of the reconnection of antiparallel vortices in the infinitely thin case and in the finite thickness case}

\author[1]{Francisco de la Hoz}
\author[2]{Sergei Iakunin\thanks{{Corresponding author: Sergei Iakunin, siakunin@bcamath.org}}}
\affil[1]{University of the Basque Country (UPV/EHU)}
\affil[2]{BCAM - Basque Center for Applied Mathematics}
\date{}

\begin{document}

    \maketitle


\begin{abstract}
    The reconnection of vortices is an important example of the transition from laminar to turbulent flow. 
    The simplest case is the reconnection of a pair of antiparallel line vortices, e.g., condensation trails of an aircraft. 
    The vortices first undergo long wave deformation (Crow waves), and then reconnect to form coherent structures.
    Although the behavior of the vortices before and after the reconnection can be clearly observed, what happens during 
    the reconnection still needs to be explained. One of the challenges is related to the fact that the vortices have finite thickness, and therefore, the 
    time and the point of the reconnection cannot be determined. Moreover, the smallest scale of coherent structures that can be 
    observed also depends on the vortex thickness. In this paper, we consider an infinitely thin vortex approximation to study
    the reconnection process.  We show that, in this case, the behavior after the reconnection is quasi-periodic, with the quasi-period being 
    independent of the angle between the vortices at the time of the reconnection. We also show that, in the Fourier 
    transform of the trajectory of the reconnection point, the frequencies that correspond to squares of integers are dominating in a similar way as in the evolution of
   a polygonal vortex under the localized induction approximation. At the end, we compare the results with a solution of the Navier-Stokes
    equations for the reconnection of a pair of antiparallel vortices with finite thickness. We use the fluid impulse to determine the 
    reconnection time, the reconnection point, and the quasi-period for this case.

\end{abstract}

\section{Introduction}

Turbulent flows are characterized by the presence of vortices that move chaotically, interact, and reconnect with each other,
forming increasingly complex structures. The basic example of this process is the reconnection of antiparallel line vortices (e.g., condensation
trails of an aircraft). The vortices initially generate a laminar flow; however, being embedded into this flow, they deform and reconnect,
giving rise to a cascade of coherent structures with rich behavior. Much of this behavior is determined by the configuration of 
the vortices at the reconnection time, which, at its turn, is not well understood. It happens due to the finite thickness of vortices,
which sets a lower limit to the scales of the coherent structures that we can observe, and to the precision with which the reconnection time
can be determined. In order to avoid this restriction, we consider infinitely thin vortices that evolve under the
localized induction approximation. We study an eye-shaped vortex that resembles the configuration close to the reconnection time,
and establish that its behavior is very reminiscent of the vortex in the shape of a regular polygon considered in~\cite{delahoz2014}. In particular,
the evolution of an eye-shaped vortex is quasi-periodic. Furthermore, we add an interaction term having the shape proposed in~\cite{klein1991a,iakunin2023},
and generate an eye-shaped vortex starting from initially smooth configuration. At the end, we compare the obtained results with a solution of 
the Navier-Stokes equations, and observe that some features, such as quasi-periodicity and vortex separation 
rate, are presented even in the case of finite thickness vortices. 

\begin{figure}
    \centering
    \includegraphics[width=0.47\textwidth]{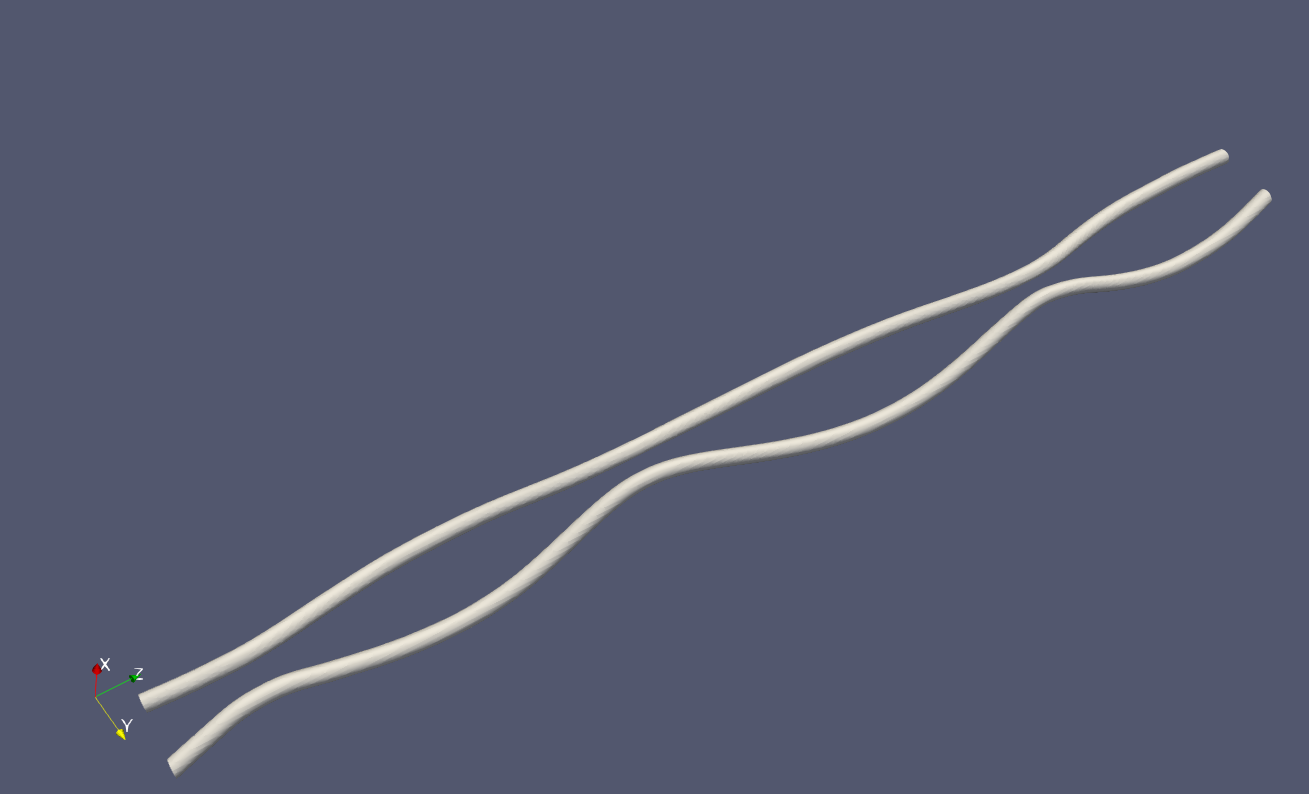}
    \includegraphics[width=0.47\textwidth]{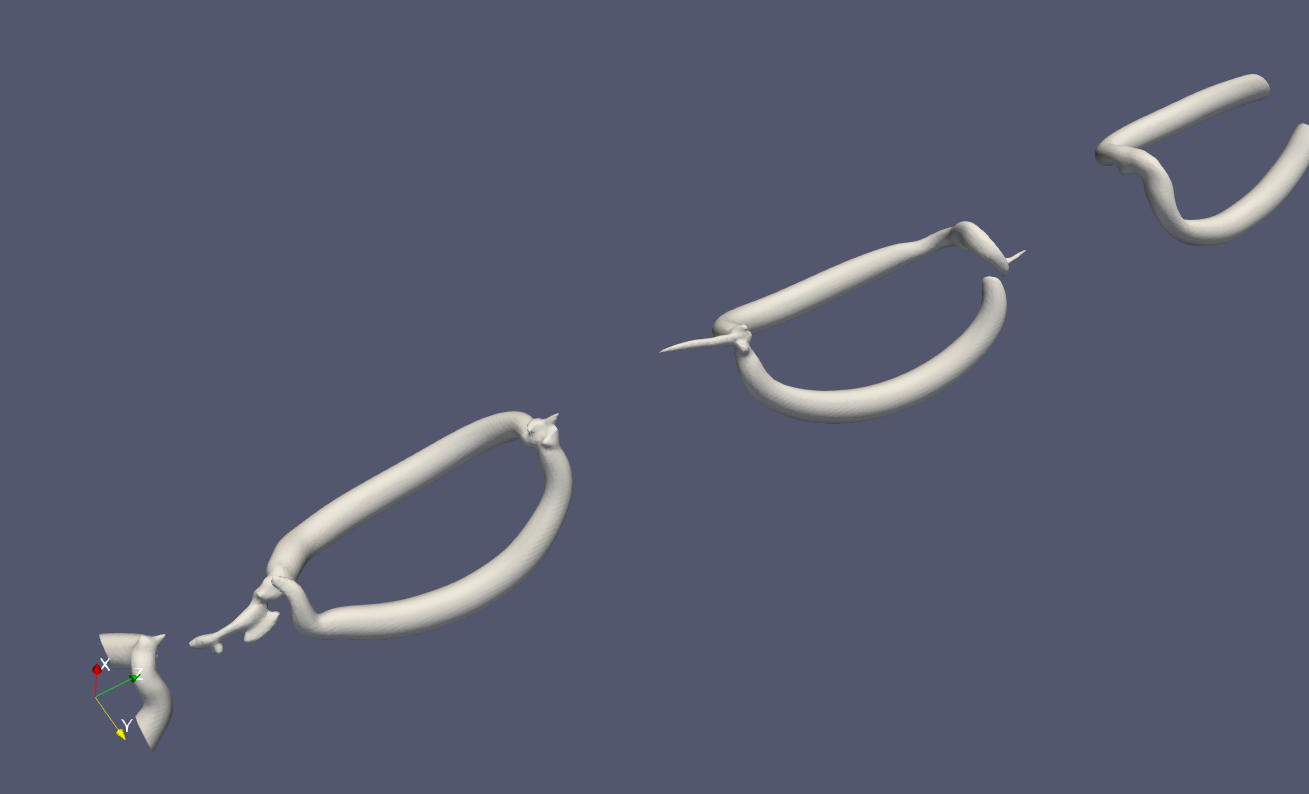}
    \caption{Pressure isosurfaces for vortex reconnection modeled in OpenFOAM.}\label{fig:ns_example}
\end{figure}

The reconnection phenomena was studied by many authors in analytical, numerical, and experimental ways. In~\cite{crow1970}, the analytical explanation of the instability of a pair of antiparallel vortices is given. If the vortices undergo a 
small perturbation, they start to rotate around their central lines with velocity proportional to the square
of the frequency of this perturbation. At the same time, the points of the vortex are pushed by the flow produced by its neighbor, and also those points that lie closer
are moving with higher velocity. In~\cite{crow1970}, it is shown that there is a frequency of perturbation when it becomes unstable and
starts to grow, leading to Crow waves. The length of these waves is proportional to the initial distance between vortices, and their orientation 
is symmetric respect to the central plane between vortices, and such that the planes containing the vortices form a corner of approximately $\pi / 4$ radians (see left side of Figure \ref{fig:ns_example}). 
It is proven that there is only long wave instability; however, in an experiment, it is possible to see also short waves~\cite{laporte2002}. These waves
are called Kelvin waves and have different nature. They are instabilities in the vortex core~\cite{ryan2007}, therefore they disappear in the approximation of infinitely
thin vortices. They also are asymmetric. The Navier-Stokes simulations~\cite{hussain2020,pontin2018} show that, after the reconnection, coherent structures emerge, which
contain a circular horseshoe and helical waves running from the reconnection point (see right side of Figure \ref{fig:ns_example}). The instability of the vortex core makes it quite challenging to analyze 
these structures, because there is no way to distinguish the deformation of the vortex central line and the deformation of the vortex core. In \cite{hussain2020}, apart from  the usage of vorticity isosurfaces and the $\lambda_2$-criterion for the visualisation of vortices, the rate of total kinetic energy dissipation is considered to define the reconnection time; it was noticed that the dissipation rate has a peak around that time.
At the end of this article, we present a method based on the analysis of the fluid impulse, which allows to study some properties of the vortex behavior after 
the reconnection, without extracting its central line.

\begin{figure}
    \centering
    \includegraphics[width=0.47\textwidth]{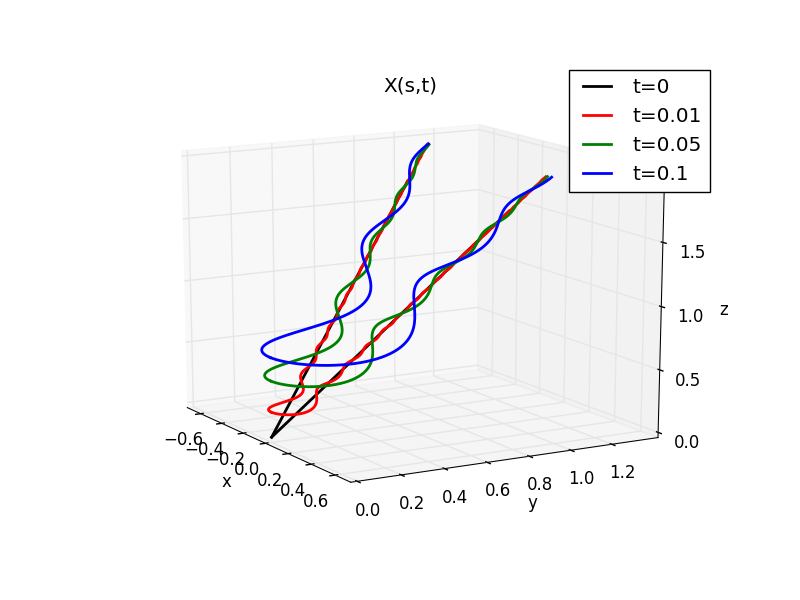}
    \caption{Evolution according to LIA of a vortex filament initially given by two half-lines that meet at a point.}\label{fig:corner}
\end{figure}

The configuration with the horseshoe and helical wave is reminiscent of ones observed in the evolution of a vortex filament that moves according 
to the localized induction approximation (LIA), and at the initial time, it is given by two half-lines that meet at point (the corner) with an angle $\theta$~\cite{vega2003}. 
This evolution, shown in Figure~\ref{fig:corner}, is self-similar, with the curvature $c_0 / \sqrt{t}$ that depends only on time $t$. The parameter $c_0$ is determined by the angle of the initial corner, according to the formula $\sin{\frac{\theta}{2}} = e^{-\pi c_0^2 / 2}$. Thus, the behavior of the corner is completely determined by its configuration at the initial moment. One can expect that the behavior after the reconnection, up to some approximation,
can be determined by the configuration at the reconnection time. However, in this article, we show that the value of the angle plays a minor role
in the subsequent behavior, the most important part being related to the symmetries of the vortex configuration. In order to prove it, we consider simpler models.

The natural example of a closed vortex that has corners is a polygonal vortex, and has been widely studied in~\cite{delahoz2014}. It is known that the behavior
of a vortex with the shape of a regular polygon with $M$ sides is periodic with period $2 \pi / M^2$. Furthermore, at times that are rational multiples of this
period, the configuration of the vortex is also a closed sequence of intervals, which is in general a skew polygon. In the Fourier transform of the corner trajectory, 
the coefficients that correspond to the squares of integers are dominating. Moreover, when $M$ tends to infinity, the trajectory is planar and tends to Riemann's non-differential function (RNDF)~\cite{banica2022}:
\begin{equation}
    \mathcal{R}(t) = \sum_{k=1}^\infty \frac{e^{i t k^2}}{k^2},\label{eq:rndf}
\end{equation}
where $i^2 = -1$. In~\cite{jaffard1996}, it was proven that RNDF is a multi-fractal. Hence, the study of this function and its relation to the behavior of vortices is quite beneficial in order to explain the multi-fractal properties of the turbulence.  This remains an open question in fluid dynamics~\cite{turiel2021}. The configuration of vortices at the reconnection has only two corners, and it is not planar, which
makes the analysis of the trajectory of the reconnection point more delicate. Therefore, we consider the eye-shaped vortex and 
its polygonal approximation. It is possible to prove that the evolution of such a vortex is quasi-periodic in the sense that, after the time $\pi / 2$, we obtain the eye-shaped vortex with the same orientation, but with different corners and different slopes. It is surprising that the quasi-period is independent of the initial corner. The analysis of the Fourier transform of the corner trajectory shows the domination of frequencies that correspond to squares of integers, similarly as in the evolution of the polygon and in RNDF. We will use this criterion
as an indicator that multi-fractal properties may be presented. In the case of the eye-shaped vortex, the corner is set as an initial condition. Thus, one may ask whether it is possible to generate such a configuration starting from a smooth one. 
The way to do it is by adding an interaction term that makes the model more complicated.

The model of interaction of almost parallel vortices that includes both local self-induction and external force is considered in a series of papers~\cite{klein1991a,klein1995}.
In a previous study~\cite{iakunin2023}, we generalized this model, and added a regularization parameter which allows us to go beyond the reconnection time, and observe coherent 
structures. However, it is impossible to follow the reconnection point for this model, because of the presence of a bridge between horseshoes. 
In this article we present a heuristic approach that allows to get rid of the bridge, using the proven fact that the behavior of the vortices after the reconnection
is almost independent from the angle that they form at the reconnection time. With this approach, we can analyze the trajectory of the reconnection point, and 
observe the domination of the frequencies that correspond to the squares. However, the shape of vortices at the reconnection time differs from the shape of the eye. For that
time, it is already possible to see a cusp instead of a corner, and the smallest horseshoe, whose size depends on the regularization parameter $r_c$. This may indicate
that the dominating squares in the Fourier transform of the reconnection point trajectory are more stable, in the sense that it is not necessary to generate a corner, but it suffices 
to provide a suitable approximation of it. For the vortices of finite thickness, the definition of the reconnection point is still unclear. Therefore, we can focus on an
integral quantity, e.g., the fluid impulse. In~\cite{iakunin2023}, it is noted that the change of fluid impulse behavior from monotone to oscillatory may be considered
as a criterion for the reconnection. It happens not only for LIA of infinitely thin vortices, but also for a solution of the Navier-Stokes equations, when the vortices
have finite thickness. Moreover, the analysis of the fluid impulse in this case allows not only to define the reconnection time and the reconnection point, but also to estimate 
the quasi-period and the vortex separation rate.

The paper has the following structure. In Section~\ref{sec:eye-shaped}, we present LIA for the evolution of an eye-shaped vortex, the most important 
features of its evolution and some conserved quantities. In the next section, i.e., Section~\ref{sec:reconnection}, we study the reconnection of infinitely thin
vortices and its relation to the eye-shaped vortex. We also introduce a heuristic approach, which allows to reconnect the vortices, get rid of the bridge, and
recover the multi-fractal behavior of the reconnection point trajectory. Later, in Section~\ref{sec:navier-stokes}, we consider the simulation of 
the reconnection of finite-thickness vortices using the Navier-Stokes equations. In this case, the fluid impulse can be used in order to determine 
the reconnection time and the reconnection point; however, there are no high frequency oscillations, unlike in the infinitely thin case. Finally, in Section~\ref{sec:conclusions}, we point out the main conclusions and indicate some further research directions.
    

\section{Localized induction approximation for the eye-shaped vortex}~\label{sec:eye-shaped}

Assume that there is an infinitely thin isolated vortex embedded into an unbounded domain of a non-viscous incompressible fluid. Moreover, at time $t$, it is
localized around a closed curve $\mathbf{X}(s,t) \in \mathbb{R}^3$, where $s$ is the (arc-length) parameter of the curve. Then, the
velocity $\mathbf{v}(\mathbf{x}, t)$ of the flow generated by this vortex at any point $\mathbf{x}$ can be calculated by using the Biot-Savart integral:
\begin{equation}
    \mathbf{v}(\mathbf{x}, t) = \frac{\Gamma}{4\pi}\int_0^L \frac{(\mathbf{x} - \mathbf{X}(q, t))\wedge \mathbf{X}_q(q, t)}
    {{|\mathbf{x} - \mathbf{X}(q, t)|}^3}dq.\label{eq:biot-savart}
\end{equation}
Here, $\Gamma$ is the vortex strength, $L$ is the length of the curve $\mathbf{X}(s,t)$, the subscripts represent partial derivatives, $\wedge$
denotes the vector product, and $|\mathbf a|$ stands for the Euclidean length of a vector $\mathbf a\in\mathbb R^3$.
The integral~\eqref{eq:biot-savart} is singular, when $\mathbf{x}$
belongs to the curve $\mathbf{X}(s,t)$. It was shown in 1906 by Da Rios that the most singular term corresponds to the rotation of particles around
the vortex central line, whereas the next term representing the deformation of the vortex due to local effects has only a logarithmic singularity.
This approximation, called the localized induction approximation (LIA), leads, after a proper time rescaling, to the vortex filament equation (VFE):
\begin{equation}
    \mathbf{X}_t(s,t) = \mathbf{X}_s(s,t) \wedge \mathbf{X}_{ss}(s,t) = \mathbf{T}(s,t) \wedge \mathbf{T}_{s}(s,t) = \kappa(s,t) \mathbf{B}(s,t),\label{eq:vfe}
\end{equation}
where $\kappa$ is the curvature, $\mathbf{T}(s,t)$ is the tangent vector, and $\mathbf{B}(s,t)$ is the binormal vector of the curve $\mathbf{X}(s,t)$.  The derivation can be found in~\cite{shaffman1992}.

There are two properties of the VFE which are quite important for us. Firstly, this equation preserves the Euclidean length of the tangent vector, i.e., $|\mathbf{T}(s,t)|~=~1$, which means that, if we parametrize the curve $\mathbf{X}(s,t)$ by its arc-length, then the parametrization will hold for all times. It implies the absence of the vortex stretching
phenomenon, i.e., all the parts of the vortex have the same strength. This assumption is not true for all cases, and, in particular, for the reconnection. Secondly, the VFE can be turned into the nolinear Schr\"odinger equation (NLSE) by the Hasimoto transform~\cite{hasimoto1972}:
\begin{equation}
    \psi(s,t) = \kappa(s,t) e^{i\int_0^s \tau(q,t) dq},\label{eq:psi_def}  
\end{equation}
where $\kappa$ and $\tau$ denote respectively the curvature and the torsion of the curve $\mathbf{X}(s,t)$. If $\mathbf{X}(s,t)$ is a solution of~\eqref{eq:vfe}, then the function
$\psi$ is a solution of NLSE:
\begin{equation}
    \psi_t(s,t) = i \psi_{ss}(s,t) + \frac{i}{2} \left({|\psi(s,t)|}^2 + A(t)\right)\psi(s,t),\label{eq:nlse}
\end{equation}
where $A(t)$ is a certain real constant that depends on time.

The equation~\eqref{eq:vfe} has been widely studied in~\cite{delahoz2014,delahoz2018} for the case when the initial condition is a regular polygon with $M$ sides. For this case,
the initial condition for $\psi(s,t)$ is just a sum of Dirac's $\delta$ functions:
\begin{equation}
    \psi_M(s,0) = \frac{2\pi}{M}\sum_{j \in \mathbb{Z}} \delta\left(s - \frac{2\pi j}{M}\right).\label{eq:poly_ic}
\end{equation}
Then, using the Galilean invariance of NLSE, it is possible to find the solution for this initial condition 
\begin{equation}
    \psi_M(s,t) = \hat{\psi}(0,t)\sum_{k \in \mathbb{Z}} e^{-i {(M k)}^2 t + i M k s},\label{eq:poly_sol}
\end{equation}
where $\hat{\psi}(0,t)$ is a constant that depends on time, and that has to be chosen in such a way that the curve $\mathbf{X}(s,t)$ is closed for all times $t$.
It is possible to show that the solution~\eqref{eq:poly_sol} is periodic in time with period $T = 2\pi / M^2$. Furthermore, at times that are rational multiples of
the time period, the solution is in general a skew polygon. However, at half a period, we have also a regular polygon with $M$ sides, but it
is rotated by an angle $\pi / M$. This section is devoted to applying the same technique to a case when the initial condition is an eye-shaped vortex.

\subsubsection*{The definition of the eye-shaped vortex}
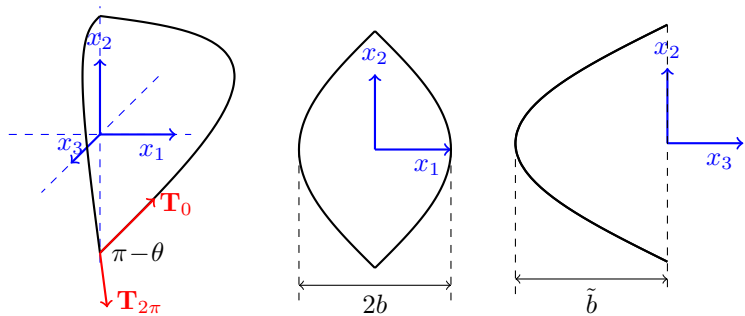
\begin{figure}
    \centering
    \begin{tikzpicture}
        \draw[blue, thick, ->] (0,0,0) -- (1,0,0) node[anchor=north east]{$x_1$}; 
        \draw[blue, thick, ->] (0,0,0) -- (0,1,0) node[anchor=south]{$x_2$};
        \draw[blue, thick, ->] (0,0,0) -- (0,0,1) node[anchor=south]{$x_3$};

        \draw[blue, dashed] (-1.2,0,0) -- (1.2,0,0); 
        \draw[blue, dashed] (0,-1.7,0) -- (0,1.7,0);
        \draw[blue, dashed] (0,0,-2) -- (0,0,2);

        \draw[domain=0:3.1415, thick, smooth, variable=\x] plot ({sin(deg(\x))}, {\x - 3.1415 / 2.}, {-2 * sin(deg(\x))});
        \draw[domain=3.1415:6.2832, thick, smooth, variable=\x] plot ({sin(deg(\x))}, {3. * 3.1415 / 2. - \x}, {2 * sin(deg(\x))});

        \draw[red, thick, ->] (0., {-3.1415 / 2.}, 0.) -- ({1. / sqrt(6)}, {-3.1415 / 2. + 1. / sqrt(6)}, {-2. / sqrt(6)}) node[pos=0.9, anchor=west] {$\mathbf T_0$};
        \draw[red, thick, ->] (0., {-3.1415 / 2.}, 0.) -- ({1.  / sqrt(6)}, {-3.1415 / 2. - 1.  / sqrt(6)}, {2. / sqrt(6)}) node[pos=0.9, anchor=west] {$\mathbf T_{2\pi}$};

        \node at ({1. / 2.}, {-3.1415 / 2.}, {0.}) {$\pi\!-\!\theta$};
    \end{tikzpicture}
    \qquad    
    \begin{tikzpicture}
        \draw[blue, thick, ->] (0,0) -- (1,0) node[anchor=north east]{$x_1$}; 
        \draw[blue, thick, ->] (0,0) -- (0,1) node[anchor=south]{$x_2$};
        
        \draw[domain=0:3.1415, thick, smooth, variable=\x] plot ({sin(deg(\x))}, {\x - 3.1415 / 2.});
        \draw[domain=3.1415:6.2832, thick, smooth, variable=\x] plot ({sin(deg(\x))}, {3. * 3.1415 / 2. - \x});        

        \draw[dashed] (-1., 0.) -- (-1., -2);
        \draw[dashed] (1., 0.) -- (1., -2);
        \draw[<->] (-1., -1.8) -- (1., -1.8) node[pos=0.5, anchor=north] {$2 b$};
    \end{tikzpicture}
    \qquad    
    \begin{tikzpicture}
        \draw[blue, thick, ->] (0,0) -- (1,0) node[anchor=north east]{$x_3$}; 
        \draw[blue, thick, ->] (0,0) -- (0,1) node[anchor=south]{$x_2$};

        \draw[domain=0:3.1415, thick, smooth, variable=\x] plot ({-2 * sin(deg(\x))},{\x - 3.1415 / 2.} );
        \draw[domain=3.1415:6.2832, thick, smooth, variable=\x] plot ({2 * sin(deg(\x))},{3. * 3.1415 / 2. - \x});
        
        \draw[dashed] (0, {3.1415 / 2.}) -- (0, -2.);
        \draw[dashed] (-2., 0.) -- (-2, -2);
        \draw[<->] (-2, -1.8) -- (0, -1.8) node[pos=0.5, anchor=north] {$\tilde b$};
    \end{tikzpicture}
    \caption{Eye-shaped vortex.}\label{fig:definition_of_eye}
\end{figure}

The eye-shaped vortices are interesting to us because they look similar to the configuration of antiparallel vortices after the reconnection, and also
because they can be considered as ``polygons'' with two sides. The initial configuration of the eye-shaped vortex is given by
\begin{equation}
    \mathbf{X}(s,0) = \begin{pmatrix}
        b \sin{s} \\ s - \pi / 2 \\ -\tilde{b}\sin{s}       
    \end{pmatrix},\ s\in(0,\pi], \quad
    \mathbf{X}(s,0) = \begin{pmatrix}
        b \sin{s} \\ 3 \pi / 2 - s \\ \tilde{b}\sin{s}       
    \end{pmatrix},\ s\in(\pi,2 \pi),\label{eq:eye_vortex}
\end{equation}
where $2b > 0$ is the maximum distance between the opposite sides of the eye, and $\tilde{b}$ is the deviation from the planar curve, see Figure~\ref{fig:definition_of_eye}. Performing the scalar product between $\mathbf T_0 = \mathbf T(0^+,0)$ and $\mathbf T_{2\pi} = \mathbf T(2\pi^-,0)$, it follows that the angle $\theta$ is given by
$$
\cos \theta = -\cos(\pi - \theta)= -\frac{\mathbf T_0 \cdot \mathbf T_{2\pi}}{|\mathbf T_0| |\mathbf T_{2\pi}|} = \frac{1 + \tilde{b}^2 - b^2}{1 + \tilde{b}^2 + b^2}.
$$
We can parametrize the curve~\eqref{eq:eye_vortex} by its arc-length and rescale it in 
such a way that the total length of the curve is equal to $2\pi$. 

\begin{figure}
    \centering
    \includegraphics[width=0.49\textwidth,height=0.392\textwidth]{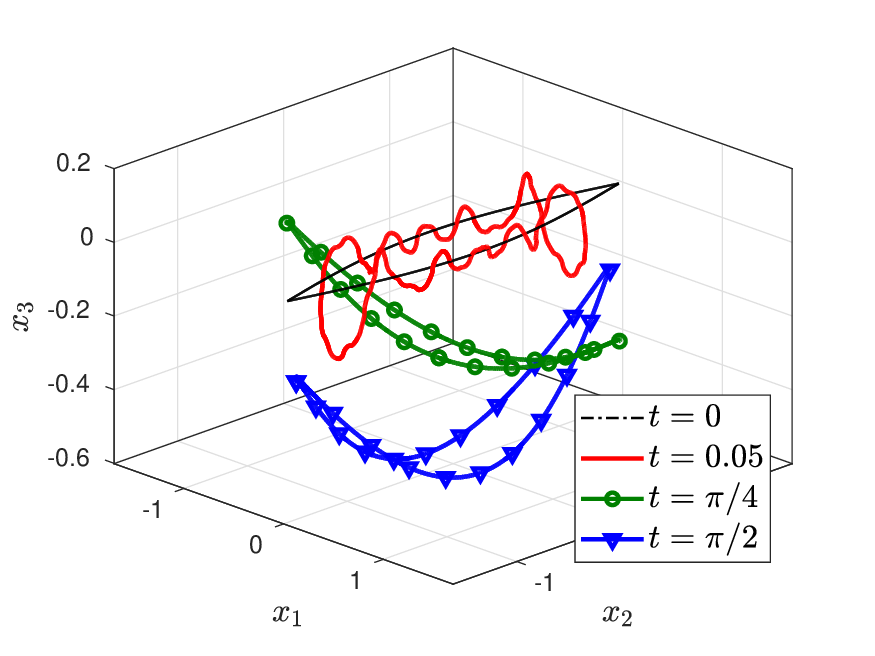}    
    \includegraphics[width=0.49\textwidth,height=0.392\textwidth]{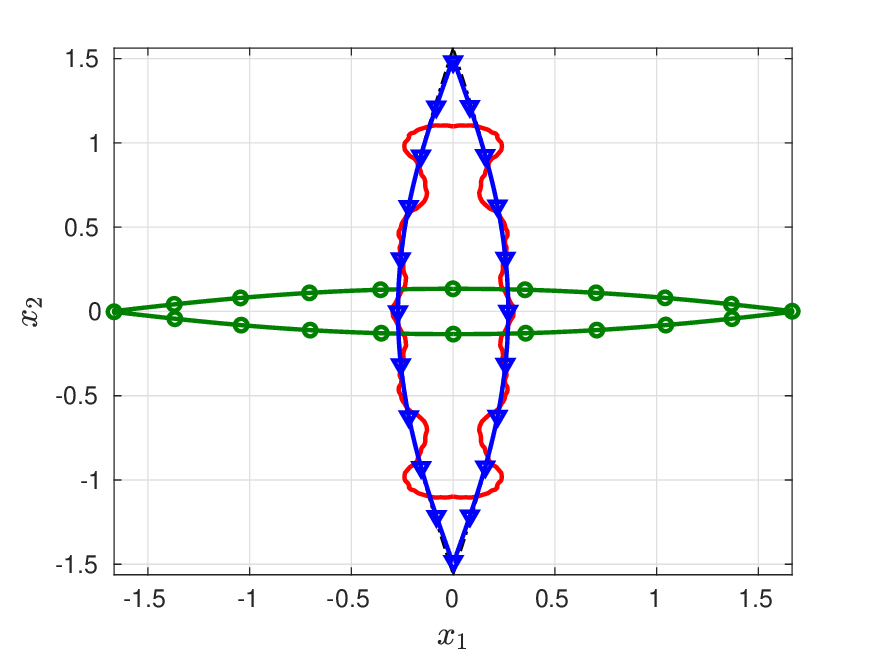}
    \caption{The evolution of the eye-shaped vortex~\eqref{eq:eye_vortex} at different time moments (initial moment, $t = 0.05$, half of a period $t = \pi / 4$ and 
    a period $t = \pi / 2$)}
    \label{fig:eye_shaped_vortex}
\end{figure}

\begin{figure}
    \centering
    \includegraphics[width=0.49\textwidth,height=0.392\textwidth]{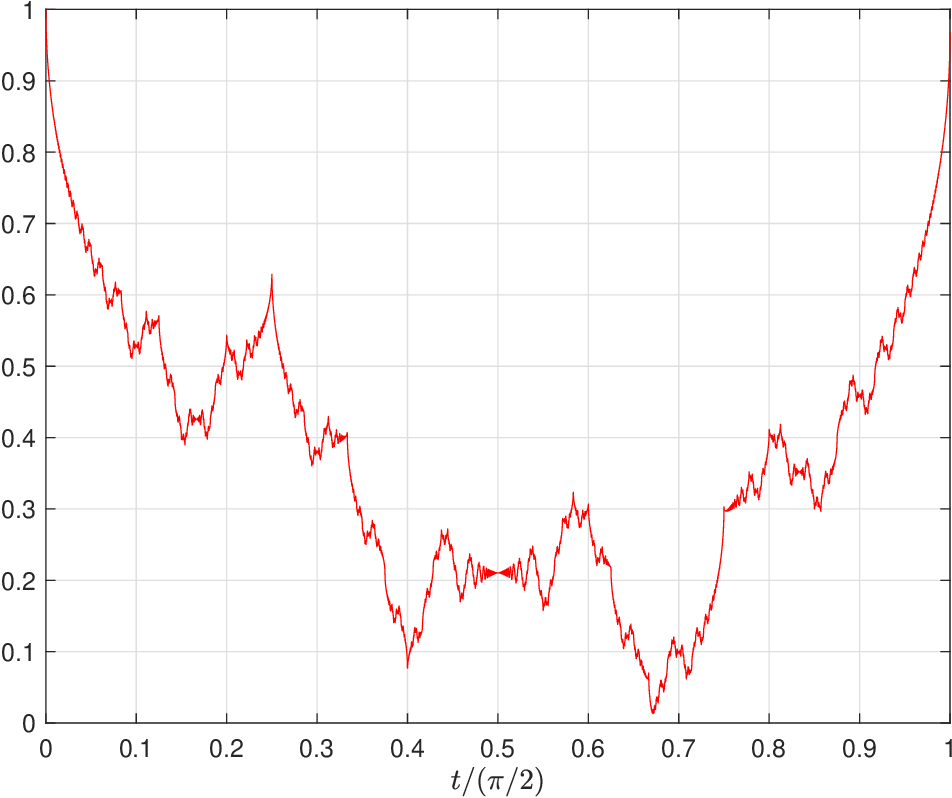}\hfill\includegraphics[width=0.49\textwidth,height=0.392\textwidth]{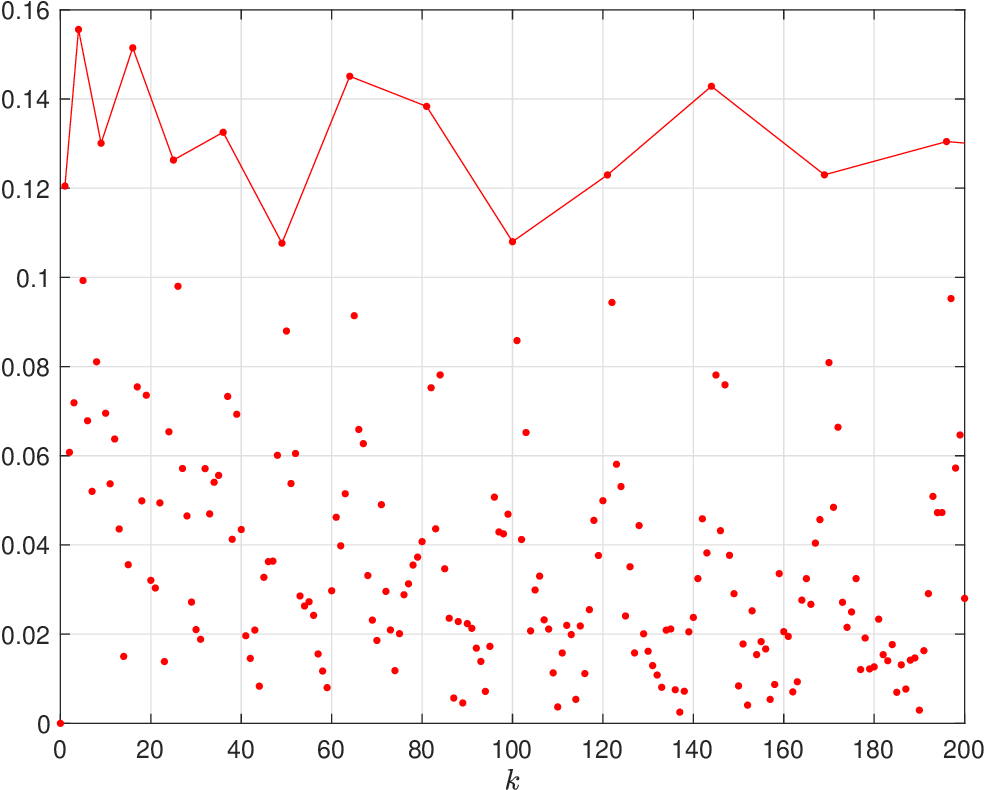}
    \caption{The modulus of the trajectory of the corner of the eye-shaped vortex $|\mathbf{X}(0,t)|$ (left) and its Fourier coefficients
    $k \widehat{|\mathbf{X}(0,t)|}(k)$ (right). The line connects the coefficients corresponding to the squares of integers.}
    \label{fig:eye_shaped_traj}
\end{figure}

\begin{figure}
    \centering
    \includegraphics[width=0.49\textwidth,height=0.392\textwidth]{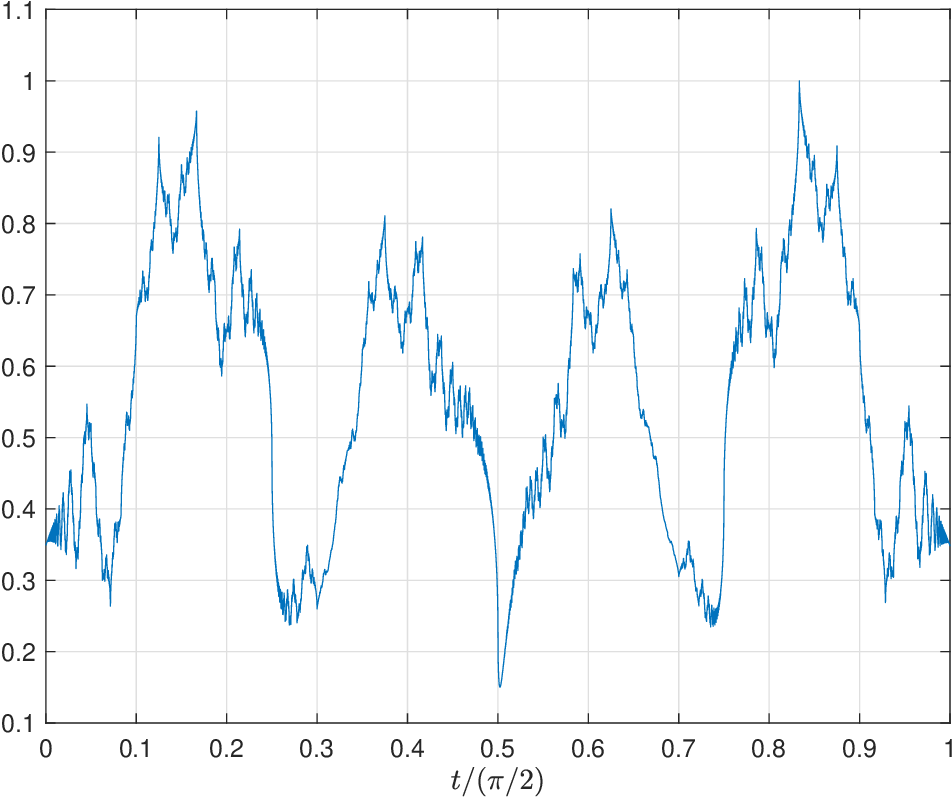}\hfill\includegraphics[width=0.49\textwidth,height=0.392\textwidth]{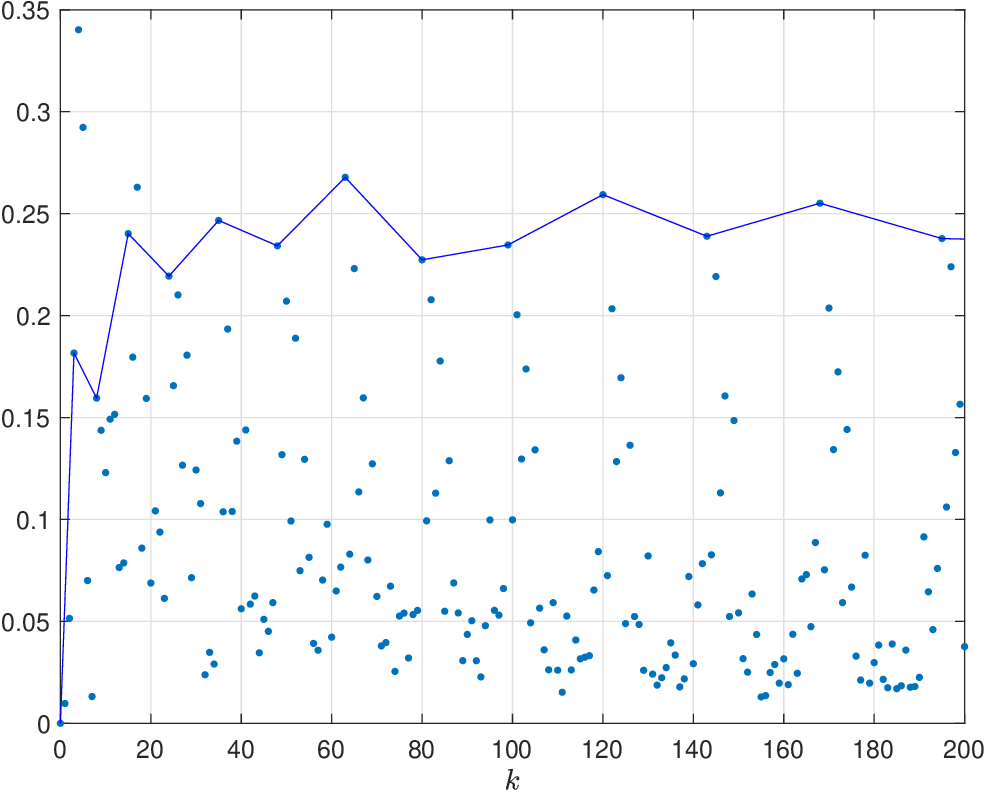}
    \caption{The modulus of the fluid impulse $|\mathbf{F}_{\pi}(t)|$ given by~\eqref{eq:fi_general_def} (left) and its Fourier coefficients
    $k \widehat{|\mathbf{F}_{\pi}(t)|}(k)$ (right). The line connects the coefficients corresponding to the squares of integers.}
    \label{fig:eye_shaped_fi}
\end{figure}

The evolution of the eye-shaped vortex~\eqref{eq:eye_vortex} is shown in Figure~\ref{fig:eye_shaped_vortex}. 
It is possible to see that the movement is quasi-periodic, since the vortex returns to the eye-shaped configuration with the same orientation (blue line with triangles), but with slightly different parameters from those at the initial time (red line), and not planar anymore. 
We can also see that, at half a period, the vortex has the shape of an eye, but rotated (green line with circles), similarly to the polygonal vortex~\cite{delahoz2018}. 
We can follow the moduli of the trajectory of the corner and of the fluid impulse $\mathbf{F}_l(t)$, where the latter is defined by:
\begin{equation}
    \mathbf{F}_l(t) = \frac{1}{2}\int_{-l/2}^{l/2} \mathbf{X}(s,t) \wedge \mathbf{T}(s,t) ds.\label{eq:fi_general_def}
\end{equation}
These two quantities are depicted in Figures~\ref{fig:eye_shaped_traj} and~\ref{fig:eye_shaped_fi}, respectively. It is possible
to see a complex structure that is reminiscent of RNDF in \eqref{eq:rndf}. It is also possible to see that the Fourier coefficients that correspond
to the squares of the integers are dominating for both the trajectory and the fluid impulse. 

\subsubsection*{The quasi-period of the eye-shaped vortex}

\begin{figure}
    \centering
    \begin{tikzpicture}
        \def\M{12} 
        \def\radius{1.5cm} 
        
        \foreach \i in {1,...,\M} {
            \draw ({\radius * sin(360. / \M * (\i - 1))}, {\radius * cos(360  / \M * (\i - 1))}) -- ({\radius * sin(360  / \M * \i)}, {\radius * cos(360  / \M * \i)});
        }
        \node[anchor = south] (1) at (0., \radius) {$1$};
        \node[anchor = south] (2) at ({\radius * sin(360  / \M)}, {\radius * cos(360  / \M)}) {$2$};
        \node[anchor = north] (M2) at (0., -\radius) {$M / 2$};
        \node[anchor = south] (M) at ({-\radius * sin(360  / \M)}, {\radius * cos(360  / \M)}) {$M$};
        
    \end{tikzpicture}
    \qquad
    \begin{tikzpicture}
        \def\M{12} 
        \def\radius{1.5cm} 
        \def\K{4}
        \foreach \i in {1,...,\M} {
            \draw[dashed] ({\radius * sin(360. / \M * (\i - 1))}, {\radius * cos(360  / \M * (\i - 1))}) -- ({\radius * sin(360  / \M * \i)}, {\radius * cos(360  / \M * \i)});
        }
        \foreach \i in {1,...,\K} {
            \draw[red, thick] ({\radius * sin(360. / \M * (\i - 1))}, {\radius * cos(360  / \M * (\i - 1))}) -- ({\radius * sin(360  / \M * \i)}, {\radius * cos(360  / \M * \i)});
            \draw[red, thick] ({\radius * sin(360. / \M * (\i + \M / 2 - 1))}, {\radius * cos(360  / \M * (\i + \M / 2  - 1))}) -- 
                                ({\radius * sin(360  / \M * (\i + \M / 2))}, {\radius * cos(360  / \M * (\i  + \M / 2)});
        }
        \node[anchor = south, red] (1) at (0., \radius) {$1$};
        \node[anchor = south, red] (2) at ({\radius * sin(360  / \M)}, {\radius * cos(360  / \M)}) {$2$};
        \node[anchor = east, red] (K1) at ({\radius * sin(360  / \M * \K)}, {\radius * cos(360  / \M * \K)}) {$K\!+\!1$};
        \node[anchor = north, red] (M2) at (0., -\radius) {$M / 2$};
        \node[anchor = west, red] (M2K1) at ({\radius * sin(360  / \M * (\M / 2 + \K))}, {\radius * cos(360  / \M * (\M / 2 + \K))}) {$M /2\!+\!K\!+\!1$};
        
    \end{tikzpicture}
    \qquad
    \begin{tikzpicture}
        \def\M{12} 
        \def\K{4}
        \def\radius{1.5cm * \M / (2 * \K)} 
        \foreach \i in {1,...,\K} {
            \draw[red, thick] ({\radius * sin(360. / \M * (\i - 1))}, {\radius * cos(360  / \M * (\i - 1))}) -- ({\radius * sin(360  / \M * \i)}, {\radius * cos(360  / \M * \i)});
            \draw[red, thick] ({\radius * (sin(360. / \M * (\i + \M / 2 - 1)) + sin(360. / \M * \K))}, {\radius * (cos(360  / \M * (\i + \M / 2  - 1)) + 1 + cos(360. / \M * \K))}) -- 
                                ({\radius * (sin(360  / \M * (\i + \M / 2)) + sin(360. / \M * \K))}, {\radius * (cos(360  / \M * (\i  + \M / 2)) + 1 + cos(360. / \M * \K))});
        }
        \node[anchor = south, red] (1) at (0., {\radius}) {$1$};
        \node[anchor = south, red] (2) at ({\radius * sin(360  / \M)}, {\radius * cos(360  / \M)}) {$2$};
        \node[anchor = west, red] (K1) at ({\radius * sin(360  / \M * \K)}, {\radius * cos(360  / \M * \K)}) {$K\!+\!1$};

        \draw ({\radius * cos(15) / 5.}, {\radius - \radius * sin(15) / 5.}) arc (-15:-105:{\radius / 5.}) node[pos=0.5, anchor=north] {$\theta_M$};

    \end{tikzpicture}
    \caption{Polygonal approximation of the eye-shaped vortex.}\label{fig:poly_eye_shaped}
\end{figure}
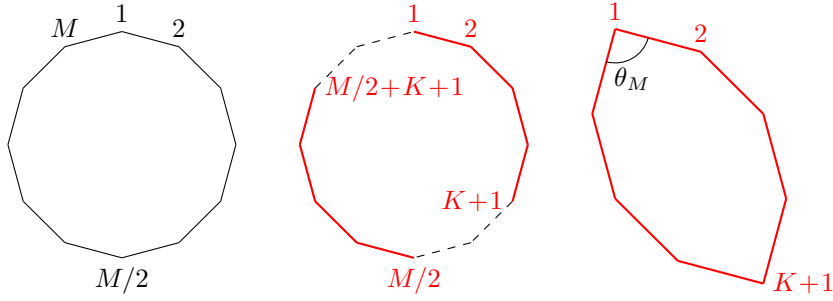

NLSE for the initial data given by~\eqref{eq:eye_vortex} does not satisfy the Galilean invariance. Therefore, the technique used for a vortex in the shape of a regular polygon is not applicable anymore. In order to find the solution, we consider a polygonal
approximation of the eye-shaped vortex. We construct this approximation in the following way (see Figure~\ref{fig:poly_eye_shaped}):
\begin{itemize}
	\item[(i)] Take a regular polygon with $M$ sides, where $M$ is even.
	
	\item[(ii)] Choose $K$ a divisor of $M$ and select the first $K$ sides from the first and second halves of the polygon.
	
	\item[(iii)] Connect the obtained curves and rescale them, so the total length equals $2 \pi$. 

\end{itemize}
Let us write $M = l K$. Then, the angles of the corners that correspond to the parts of the polygon are $2 \pi (M - 2) / M$, and the angles of the two corners that correspond to the connection between the two parts are $\theta_M = 2 \pi (M - l) / (l M)$. Note that, when $M$ tends to infinity, the parts of the polygon turn into smooth lines, whereas the corners of the eye have a limit $\theta = 2 \pi / l$. Here, we are supposing that $l$ is a
natural number; however, this construction can be done for any rational $l = M / K > 1$. In the case of $l = 2$, the eye turns into a regular polygon. 

For the vortex considered, the initial condition for NLSE~\eqref{eq:nlse} reads:
\begin{equation}
    \psi(s, 0) = \frac{\pi (l - 2)}{l} \sum_{j \in \mathbb{Z}} \delta(s - \pi j) + \frac{2 \pi}{M} \sum_{j \in \mathbb{Z}} \delta\left(s - \frac{\pi j}{K}\right).\label{eq:init_cond}
\end{equation}
Using Poisson's summation formula 
\begin{equation}
    \sum_{j\in\mathbb{Z}} f(j) = \sum_{j\in\mathbb{Z}} \hat{f}(2 \pi j) = \sum_{j\in\mathbb{Z}} \int_{-\infty}^\infty e^{-2 \pi i j y} f(y) dy,\label{eq:poisson_formula}
\end{equation}
we can rewrite the initial condition~\eqref{eq:init_cond} as
\begin{equation}
    \psi(s, 0) = \frac{(l - 2)}{l} \sum_{j \in \mathbb{Z}} e^{-2i j s} + \frac{2}{M} \sum_{j \in \mathbb{Z}} e^{-2i K j s}.\label{eq:init_cond_2}
\end{equation}
It is possible to see that the initial condition~\eqref{eq:init_cond_2} satisfies the Galilean invariance, i.e., $\psi(s, 0) = e^{-2 i K q s}\psi(s, 0)$, for all $q \in \mathbb{Z}$. This property can be used to find the solution of NLSE, similarly as it was done for a regular polygon in~\cite{delahoz2014}. Note that, if $\psi(s,t)$ is a solution of~\eqref{eq:nlse}, then, $e^{i z s - i z^2 t}\psi(s - 2 z t,t)$ is also a solution, for all $z \in \mathbb{R}$. We can choose $z = K q$ and find the coefficients of the Fourier transform of $\psi(s,t)$ with respect to $s$:
\begin{align}
    \hat{\psi}(j,t) & = \frac{1}{\pi} \int_0^\pi e^{-2 i j s} \psi(s,t) dsr
    \cr
    & = \frac{1}{\pi} \int_0^\pi e^{-2 i j s} e^{2 i K q s - 4 i (K q)^2 t}\psi(s - 4 K q t,t) ds
    \cr
    & = \frac{1}{\pi} \int_0^\pi e^{-2 i (j - Kq) (s - 4 K q) - 8 i j K q t + 8 i (Kq)^2 t}\psi(s - 4 K q t,t)
    \cr
    & = e^{4 i (K q)^2 t - 8 i j K Q t}\hat{\psi}(j - K q,t),\label{eq:fourier_transform}
\end{align}
for all $q\in\mathbb Z$. Thus, choosing $j = K q + r$, where $r \in\{ 0,\dots,K-1\}$ we obtain that there are only $K$ independent coefficients, so the solution can be written as 
\begin{equation}
    \psi(s,t) = \sum_{r = 0}^{K - 1} \hat{\psi}(r, t) e^{4 i r^2 t} \sum_{q\in\mathbb{Z}} e^{-4 i (Kq + r)^2 t + 2 i (Kq + r) s} = \sum_{r = 0}^{K - 1} \phi_r(t) \theta_r(s,t).\label{eq:solution}
\end{equation}
Note that all $\theta_r(s,t)$ are solutions of the linear Schr\"odinger equation. If we choose $t = \pi / 2$, all the exponents in \eqref{eq:solution} that depend on time become equal to $1$, so the solution is
\begin{equation}
    \psi(s,\pi / 2) = \sum_{r = 0}^{K - 1} \hat{\psi}(r, \pi / 2) \sum_{q\in\mathbb{Z}} e^{2 i (Kq + r) s} = 
    \sum_{r = 0}^{K - 1} \pi \hat{\psi}(r, \pi / 2) \sum_{q\in\mathbb{Z}} \delta\left(s - \frac{\pi r}{K} + \pi q\right).
\end{equation}
This corresponds to a polygon (not necessary planar) with $2K$ corners, similar to the one that we initially had. Furthermore, the numerical simulation suggests that two of the corners at time $\pi / 2$ that are opposite to each other are much sharper than the rest, so the shape is close to an eye. The orientation of the eye is also similar to the initial one; however, the values of the angles are different, and the curve is not planar. Note that the quasi-period  $\pi / 2$ is independent of the number of sides in the polygonal approximation and also of the angle $\theta$. This fact is important for the reconnection of vortices, because it is very involved to find the angle of the corner made by vortices. However, since the subsequent evolution of a vortex depends mostly on the number of corners and symmetries, but not on the angles, we can perform an approximate reconnection without loss of accuracy. 

\subsubsection*{Conserved quantities}
Let us consider a particular case of polygonal approximation of the eye-shaped vortex when $K = 2$. It corresponds to a rhombus, and according to \eqref{eq:solution}, the whole evolution can be described by the complex-valued functions $\hat{\psi}(0, t)$ and $\hat{\psi}(1, t)$. It is possible 
to show that, at time moments $t_k = k \pi / 8$, the solution is also a rhombus, but not necessary planar. The solution~\eqref{eq:solution} does not imply that
the curve is closed; therefore, enforcing this condition, we can find constrains for the functions $\hat{\psi}(0, t)$ and $\hat{\psi}(1, t)$.

At time moments when the solution is a rhombus, the $2\pi$-periodic tangent vector $\mathbf{T}(s,t)$ has four jumps; more precisely, in $s\in(0,2\pi)$, we have
\begin{equation}
    \mathbf{T}(s) = \left\{
    \begin{aligned}
        &\mathbf{T}_0, & & s\in(0,\pi/2),\\
        &\mathbf{T}_1, & & s\in(\pi/2,\pi),\\
        &\mathbf{T}_2, & & s\in(\pi, 3\pi/2),\\
        &\mathbf{T}_3, & & s\in(3\pi/2, 2\pi).
    \end{aligned}
    \right.
\end{equation}
Moreover, the transition between different sides is determined by the rotation matrices
\begin{equation*}
	\mathbf{M}_k = \begin{pmatrix}
		\cos{\rho_k} & \sin{\rho_k} \cos{\theta_k} & \sin{\rho_k} \sin{\theta_k}\\
		-\sin{\rho_k} \cos{\theta_k} & \cos{\rho_k} \cos^2{\theta_k}\!+\!\sin^2{\theta_k} & (\cos{\rho_k}\!-\!1) \cos{\theta_k} \sin{\theta_k} \\
		-\sin{\rho_k} \sin{\theta_k} & (\cos{\rho_k}\!-\!1) \cos{\theta_k} \sin{\theta_k} & \cos{\rho_k} \sin^2{\theta_k}\!+\!\cos^2{\theta_k}
	\end{pmatrix},
\end{equation*}
where $\rho_k e^{i \theta_k}$ is determined from the coefficients that multiply the $\delta$-functions in the solution of NLSE \cite{delahoz2014}, i.e., from $\hat{\psi}(k, t)$. Then,
\begin{equation*}
	\begin{pmatrix} \mathbf{T}_{k+1} \\ \mathbf{e}_{1,k+1} \\ \mathbf{e}_{2,k+1} \end{pmatrix} =
	\mathbf{M}_k \begin{pmatrix} \mathbf{T}_{k} \\ \mathbf{e}_{1,k} \\ \mathbf{e}_{2,k} \end{pmatrix},    
\end{equation*}
where $\mathbf{e}_{1}$ and $\mathbf{e}_{2}$ generalize the normal and binormal vectors in the Frenet-Serret equations, and form an orthogonal basis of $\mathbb R^3$ with the tangent vector $\mathbf T$. On the other hand, due to the symmetries of the rhombus, only two rotation matrices, $\mathbf M_0$ and $\mathbf M_1$, need to be determined, because $\mathbf M_{2l} = \mathbf M_0$ and $\mathbf M_{2l+1} = \mathbf M_1$, for all $l$. Furthermore, their product $\mathbf M_1\mathbf M_0$ induces a rotation of $\pi$ radians, so the trace is $1 + 2\cos(\pi) = -1$, from which we conclude, after simplification, that
\begin{equation}
    \cos(\theta_1 - \theta_2) = \cot(\rho_1 / 2) \cot(\rho_2 / 2).\label{eq:trace_conservation}
\end{equation}
The derivation is presented in appendix~\ref{sec:first_cons_quant}.

Another conserved quantity is the full fluid impulse, which is defined using \eqref{eq:fi_general_def} with $l = 2\pi$:
\begin{equation}
    \mathbf{F}(t) = \frac{1}{2}\int_{-\pi}^{\pi} \mathbf{X}(s,t) \wedge \mathbf{T}(s,t) ds.\label{eq:fluid_impulse}
\end{equation}
It is easy to see that this quantity does not depend on time for the solution of VFE~\eqref{eq:vfe}. Indeed,
\begin{align*}
    \mathbf{F}_t(t) & = \frac{1}{2}\int_0^{2\pi} \mathbf{X}_t(s,t) \wedge \mathbf{T}(s,t) ds + \frac{1}{2}\int_0^{2\pi} \mathbf{X}(s,t) \wedge \mathbf{T}_t(s,t) ds 
    \cr
    & = \int_0^{2\pi} \mathbf{X}_t(s,t) \wedge \mathbf{T}(s,t) ds + \left.\mathbf{X}(s,t) \wedge \mathbf{X}_t(s,t)\right|_{s = 0}^{2\pi}
    \cr
    & = \int_0^{2\pi} \left(\mathbf{T}(s,t) \wedge \mathbf{T}_s(s,t)\right)\wedge \mathbf{T}(s,t) ds = 0.
\end{align*}
In the case of a polygon, the fluid impulse is
\begin{align*}
\mathbf{F} & = \frac12\sum_{k = 0}^{M-1}\int_{2\pi k/M}^{2\pi (k+1)/M}\mathbf X(s)\wedge\mathbf T(s)ds
\cr
& =
\frac12\sum_{k = 0}^{M-1}\int_{2\pi k/M}^{2\pi (k+1)/M}\left[\mathbf X_k + \left(s - \frac{2\pi k}{M}\right)\mathbf T_k\right]\wedge\mathbf T_kds,
\end{align*}
where $\mathbf{X}_k$ is the position of the $k$th vertex and $\mathbf{T}_k$ is the direction vector of the $k$th side. Then,
\begin{equation*}
	\mathbf{F} = \frac12\sum_{k = 0}^{M-1}\int_{2\pi k/M}^{2\pi (k+1)/M}\mathbf X_k\wedge\mathbf T_kds = \frac{\pi}{M}\sum_{k = 0}^{M-1}\mathbf X_k\wedge\mathbf T_k.
\end{equation*}
Taking into account that $\mathbf{X}_k = \mathbf{X}_{k-1} + (2\pi/M)\mathbf{T}_{k-1}$, it follows that
$$
\mathbf X_k = \mathbf X_0 + \frac{2\pi}{M}\sum_{l = 0}^{k-1}\mathbf T_l, \quad k > 0.
$$
Hence,
\begin{equation*}
	\mathbf{F} = \frac{\pi}{M}\mathbf X_0\wedge\sum_{k=0}^{M-1}\mathbf T_k + 2 \left(\frac{\pi}{M}\right)^2\sum_{k = 1}^{M-1}\sum_{l = 0}^{k-1}\mathbf T_l\wedge\mathbf T_k = 2 \left(\frac{\pi}{M}\right)^2\sum_{k = 1}^{M-1}\sum_{l = 0}^{k-1}\mathbf T_l\wedge\mathbf T_k,
\end{equation*}
where we have used that the curve is closed, i.e., $\sum_{k = 0}^{M - 1} \mathbf{T}_k = \mathbf{0}$. In the case of a rhombus with $M = 4$, the last formula becomes
\begin{align*}
\mathbf{F} & = \frac{\pi^2}{8}\sum_{k = 1}^{3}\sum_{l = 0}^{k-1}\mathbf T_l\wedge\mathbf T_k
	\cr
& = \frac{\pi^2}{8}\left(\mathbf T_0\wedge\mathbf T_1+ (\mathbf T_0 + \mathbf T_1)\wedge\mathbf T_2 + (\mathbf T_0 + \mathbf T_1 + \mathbf T_2)\wedge\mathbf T_3\right)
	\cr
& = \frac{\pi^2}{8}\left(\mathbf T_0\wedge\mathbf T_1 - (\mathbf T_2 + \mathbf T_3) \wedge\mathbf T_2 - \mathbf T_3\wedge \mathbf T_3\right)
	\cr
& = \frac{\pi^2}{8}\left(\mathbf T_0\wedge\mathbf T_1 - \mathbf T_3 \wedge\mathbf T_2\right) = \frac{\pi^2}{8}\left(\mathbf T_0\wedge\mathbf T_1 + \mathbf T_2 \wedge\mathbf T_3\right),
\end{align*}
where we have used that $\mathbf T_0 + \mathbf T_1 + \mathbf T_2 + \mathbf T_3 = \mathbf 0$. After simplification, the square of the modulus of $\mathbf F$ becomes
\begin{equation}
|\mathbf{F}|^2 = 
\frac{\pi^4}{16}\sin^2(\rho_0)\left(1 - \sin^2(\rho_1/2)\sin^2(\theta_0 - \theta_1)\right),\label{eq:fi_expression}
\end{equation}
which is, as said above, a preserved quantity. Introducing \eqref{eq:trace_conservation} into this last expression, we conclude that
\begin{equation}
(1 + \cos(\rho_0)) (1 + \cos(\rho_1)) = \const(t) \Longleftrightarrow \cos(\rho_0/2)\cos(\rho_1/2) = \const(t).\label{eq:fi_conservation}
\end{equation}
The derivation can be found in appendix~\ref{sec:second_cons_quant}.
This property has been observed numerically for general polygons \cite{delahoz2018}; however, there is still no proof.    

\section{Reconnection of infinitely thin vortices}~\label{sec:reconnection}

In the previous section, the eye-shaped vortex has been parametrized by arc-length, and the vortex stretching has been neglected. However, if different parts of the vortex have different strength, we cannot use the symmetry anymore, and the
period depends now on how the vortex is stretched. The situation is even more complicated if the non-local effects and interaction between
vortices are included, since they oblige the vortex strength to change in time. In this section, we consider the model of vortex 
reconnection introduced in~\cite{iakunin2023}. We propose a heuristic approach on how to change the boundary conditions when the vortices are close to each other, obtaining in this way a vortex shape very close to the eye considered in ~\eqref{eq:eye_vortex}. We study how the presence of the interaction term and the arc-length parametrization affect the behavior of the vortices.

\subsubsection*{Model for the interaction of a pair of antiparallel vortices}

\begin{figure}
    \centering
    \begin{tikzpicture}
        \draw[black, dashed] (0, -2.2) -- (0, 2.2);
        \draw[black, <->] (0, -1.7) -- (1.5, -1.7) node[pos=0.5, anchor=north] {$b$};
        \draw[black, ->] (0, 0) -- (1, 0) node[pos = 1, anchor=west] {$\mathbf{e}_1$};
        \draw[black, ->] (0, 0.) -- (0, 1) node[pos = 1, anchor=west] {$\mathbf{e}_3$};
        \draw[black, ->] (0, 0) -- (-0.7, -0.7) node[pos = 1, anchor=west] {$\mathbf{e}_2$};
        \draw[line width = 0.5mm, blue] (1.5, -2) -- (1.5, 2) node[pos = 1, anchor=west] {$\mathbf{X}(s,t)$};
        \draw[dashed, line width = 0.5mm, red] (-1.5, -2) -- (-1.5, 2) node[pos = 1, anchor=east] {$\mathbf{Y}(s,t)$};;
        \draw[line width = 0.5mm, blue, ->] (1.9, 1) arc(0:180:0.4);
        \draw[dashed, line width = 0.5mm, red, <-] (-1.1, 1) arc(0:180:0.4);
    \end{tikzpicture}    
    \caption{Initial configuration of vortices}\label{fig:vortex_init}
\end{figure}
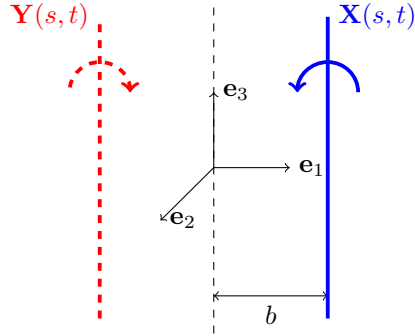

Consider a pair of infinitely thin vortices embedded into a non-viscous fluid and located initially along infinite lines going in
vertical direction, as shown in Figure~\ref{fig:vortex_init}. Suppose also that these vortices have equal strength, but rotate in opposite
directions. In the previous work~\cite{iakunin2023}, we derived a variation of VFE that models this pair of vortices:
\begin{equation}
    \mathbf{X}_t(s,t) = \frac{\mathbf{X}_s(s,t) \wedge \mathbf{X}_{ss}(s,t)}{|\mathbf{X}_s(s,t)|^3} - 
    \frac{\varepsilon x_1(s,t)}{x_1^2(s,t) + r_c^2} \frac{\mathbf{X}_s(s,t) \wedge \mathbf{e}_1}{|\mathbf{X}_s(s,t)|}, \label{eq:vfe_interaction}
\end{equation}
considering periodic boundary conditions for all components of $\mathbf{X}$. Here, $\varepsilon$ represents the strength of the interaction with 
respect to the local self-induction, and $r_c$ is a regularization parameter related to the size of the vortex core and viscosity. According to~\cite{shaffman1992}, the radius of the vortex core satisfies $r_c \sim \sqrt{\nu t}$, where $\nu$ is the viscosity and $t$ is the time.
In this model, the modulus of the tangent vector is not a constant, but depends on the distance between vortices:
\begin{equation}
    |\mathbf{X}_s(s,t)| = c(s) {(x_1^2(s,t) + r_c^2)}^{-\varepsilon / 2}, \label{eq:arclength}
\end{equation}
where $0 < c_0 \le c(s) \le C_0$ is a smooth function that depends only on $s$; therefore, it can be derived from the initial conditions. This function
varies typically very little from a constant. Thus, the interaction term leads to vortex stretching, which implies that regions of the first vortex that are closer to 
the second one have bigger circulation.

\begin{figure}
    \centering
    \includegraphics[width=0.49\textwidth]{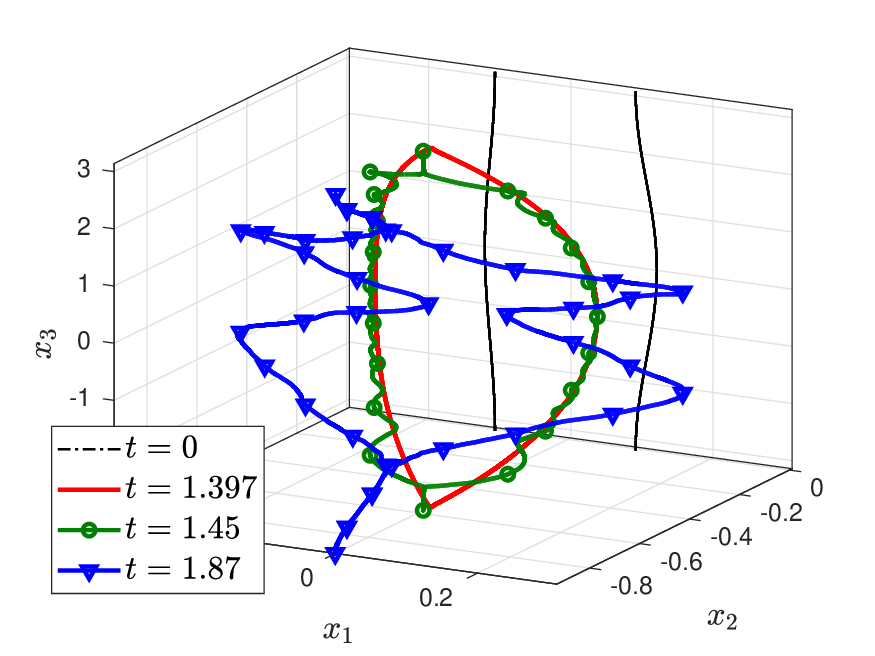}
    \includegraphics[width=0.49\textwidth]{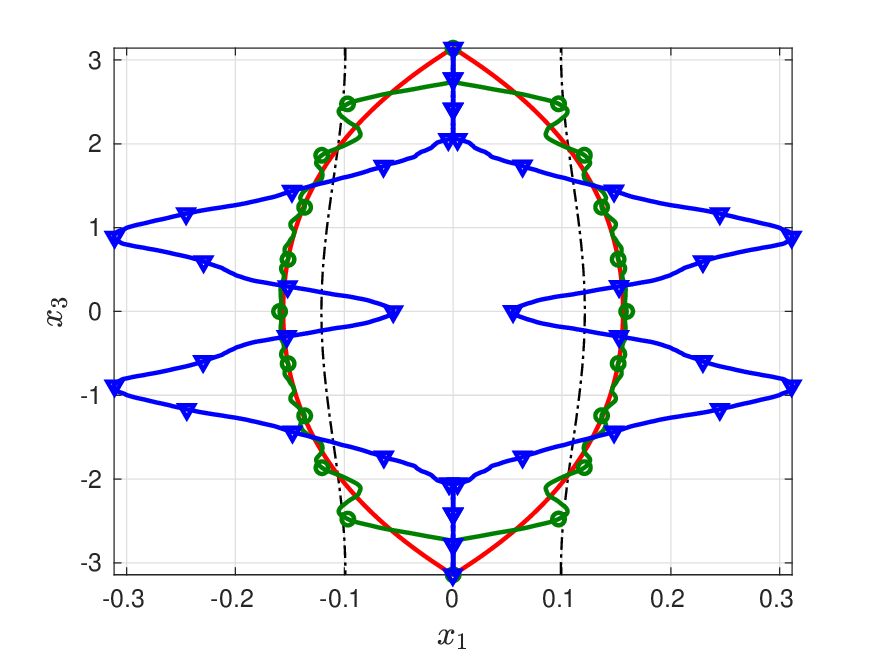}
    \caption{Configuration of the vortices at different times, for $r_c = 7.5\times 10^{-4}$, $\varepsilon = 0.03$. The initial distance is $b = 0.22$, and
    the discretization has been done taking $6000$ nodes.}
    \label{fig:crow_sol}
\end{figure}

The solution of~\eqref{eq:vfe_interaction} for $\varepsilon = 0.03$ and $r_c = 7.5\times 10^{-4}$ is depicted in Figure~\ref{fig:crow_sol} at
four time moments. In time $t = 0$, the vortices are almost parallel, and separated by a distance $2 b = 0.44$. At time $t = 1.397$, the reconnection happens, and the vortices
form a shape very close to an eye. Later on, at $t = 1.45$ a horseshoe and helical waves emerge. However, it is also possible to see a bridge that is due to
periodic boundary conditions set for each vortex. At time $t = 1.87$, it is possible to observe that the vortex has more corners and the bridge is much larger.
Eventually, the bridge produces noise, the periodicity is destroyed, and the coherent structures that appear in the evolution of the eye-shaped vortex are now absent. Furthermore,
the presence of the bridge makes it impossible to define the reconnection point and follow its trajectory. In the previous work~\cite{iakunin2023}, the fluid
impulse was used for the analysis of the vortex behavior after the reconnection. Since this quantity is an integral, it is more stable. However, the definition of the period
and the extraction of frequencies corresponding to squares of integers are still impossible in this case. In this article, we present another approach based on developing a reconnection criterion and changing the boundary conditions, when the vortices satisfy that criterion. 

\subsubsection*{Reconnection model}

We can avoid the emergence of the bridge if we change the boundary conditions when the vortices are close to each other. The main problem here consists in estimating the angle of the corner that the vortices are forming. This corner is not known, because of the regularization parameter $r_c$ that represents the smallest scales of structures that
can be observed close to the reconnection time. In Section~\ref{sec:eye-shaped}, we have proved that the quasi-period of an eye-shaped 
vortex is independent of the angle of the corner. Thus, we can perform the reconnection without an estimation of the angle, and expect that the behavior of the vortex does not change a lot.

The natural way to find the reconnection time is by means of the analysis of the distance between vortices. We can introduce a threshold 
$\text{th}_{x_1} = 10^{-6}$ and perform the reconnection at the smallest time $t_{rec}$, such that 
\begin{equation}
    \min_{s \in (0, 2\pi)} x_1(s, t_{rec}) \le \text{th}_{x_1}.\label{eq:min_u_crit}
\end{equation}
At this moment: 
\begin{enumerate}
    \item We change the boundary conditions for $x_1(s,t)$ and $x_3(s,t)$ to asymmetric:
    $$
    x_1(0,t) = -x_1(2\pi, t),\quad x_3(0,t) = -x_3(2\pi, t).
    $$
    \item We set $\varepsilon = 0$, in order to avoid a second reconnection.
    \item We reparametrize the vortices by arc-length, since it is now preserved.
\end{enumerate}

\begin{figure}
    \centering
    \includegraphics[width=0.49 \textwidth]{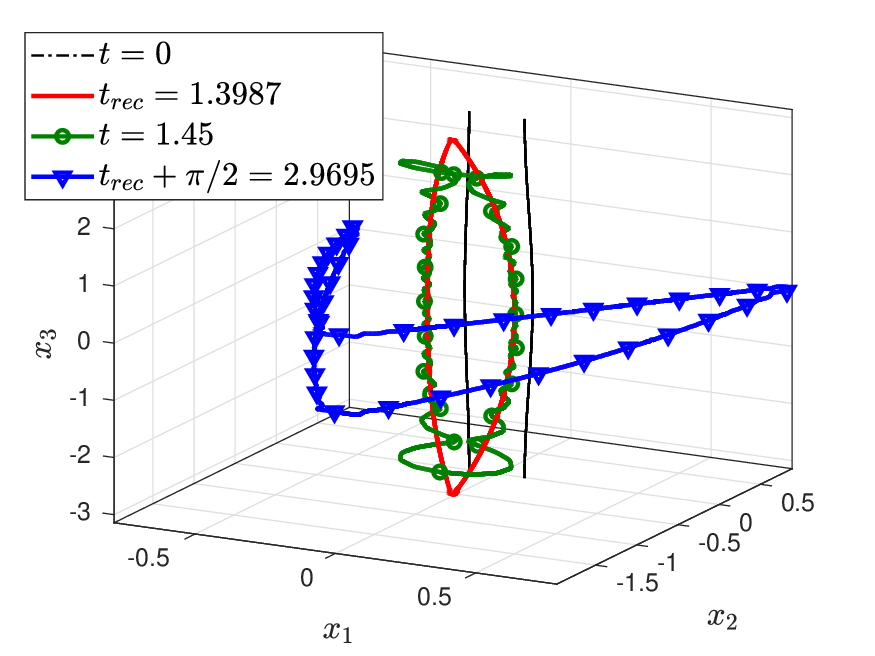}        
    \includegraphics[width=0.49 \textwidth]{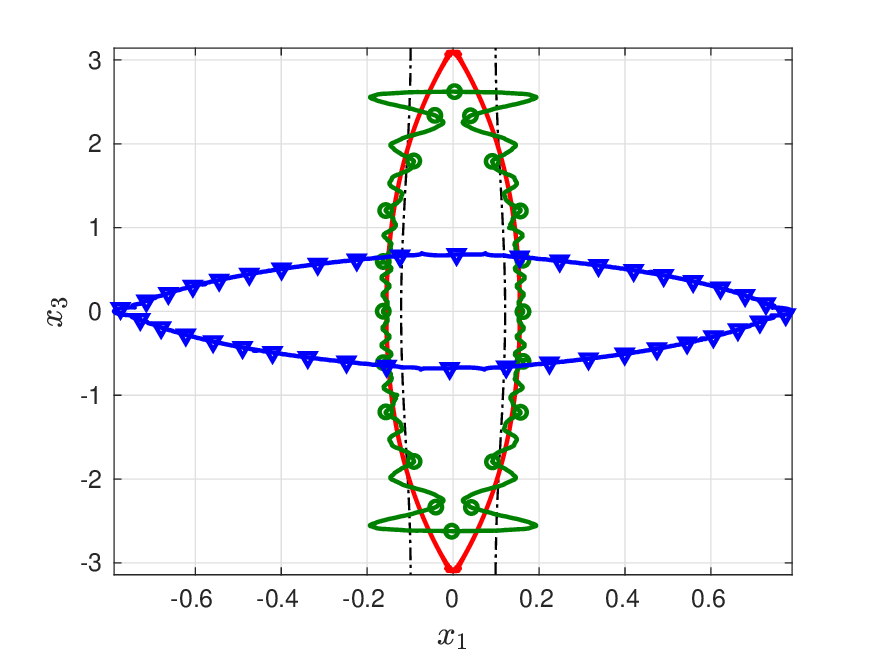}        
    \caption{Configuration of the vortices at different times, for $r_c = 7.5\times 10^{-4}$, $\varepsilon = 0.03$. The initial distance is $b = 0.22$, and the discretization has been done taking $6000$ nodes. The reconnection is allowed using the criterion~\eqref{eq:min_u_crit}, with $\text{th}_{x_1} = 10^{-6}$. 
    The interaction term is neglected, and the vortices are reparametrized by arc-length after the reconnection.} \label{fig:rec_sol}
\end{figure}

\begin{figure}
    \centering
    \includegraphics[width=0.49 \textwidth]{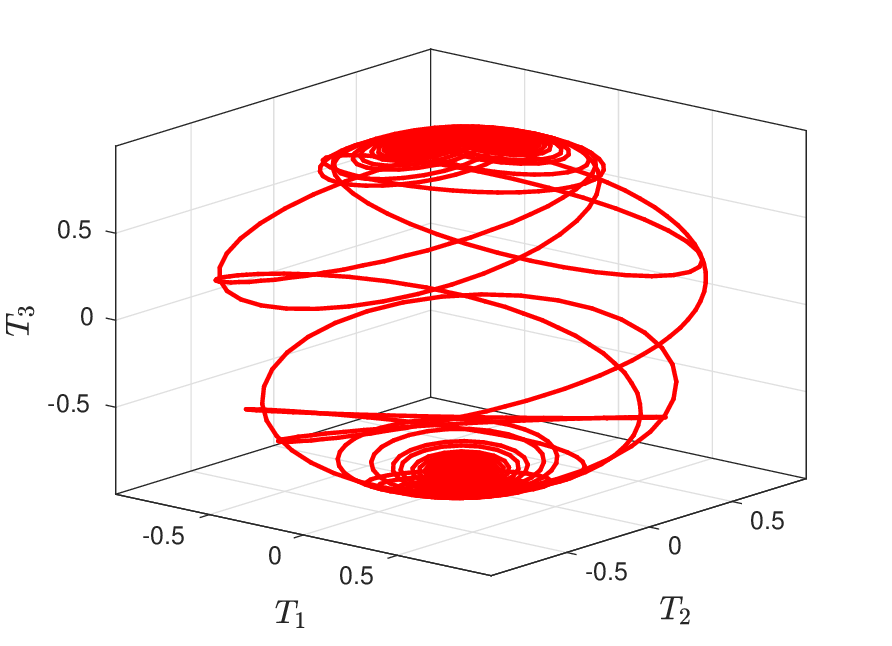}        
    \includegraphics[width=0.49 \textwidth]{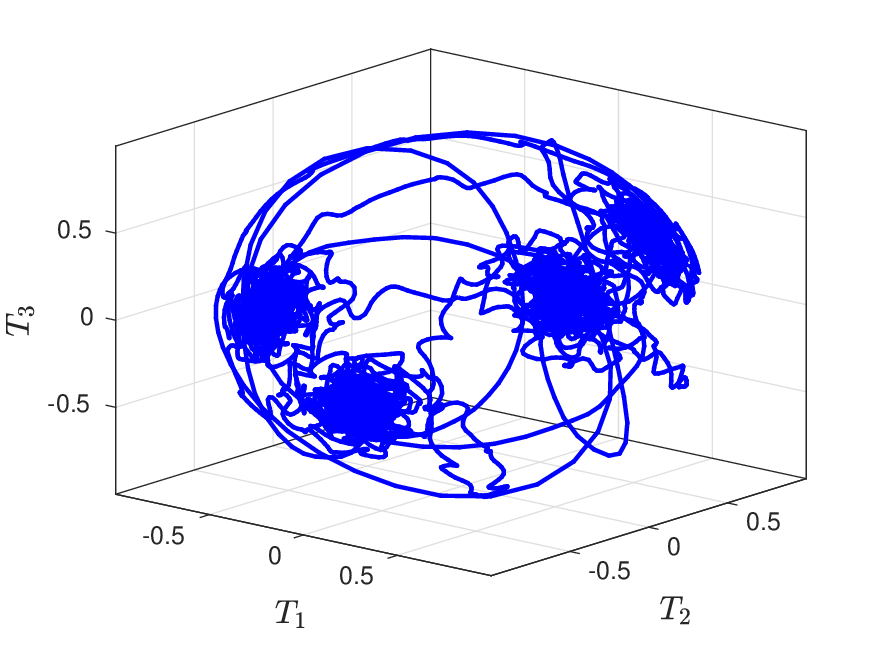}        
    \caption{Tangent vectors of the vortices depicted in Figure~\ref{fig:rec_sol} around time moments $t_{rec}$ (left) and $t_{rec} + \pi / 2$ (right). 
    The spirals are helical waves that run from the horseshoes that emerge instead of corners. The centers of the spirals correspond to the directions of the edges of the corners.}  
    \label{fig:rec_sol_tan}
\end{figure}

The evolution of the vortices with reconnection according to criterion~\eqref{eq:min_u_crit} is depicted in Figure~\ref{fig:rec_sol}.
It is possible to see that, after the reconnection, the evolution of the vortices is very reminiscent of the evolution of the eye-shaped vortex~\eqref{eq:eye_vortex}.
In particular, at time $\pi / 2$ after the reconnection, the vortex also has the shape of an eye (two corners connected by arcs), as it is depicted in Figure~\ref{fig:rec_sol_tan}.
Thus, we can expect that the behavior is quasi-periodic with period $\pi$. This increment of period with respect to the eye-shaped case is due to the
increment of the vortex length: $4\pi$ instead of $2 \pi$. 

\begin{figure}
    \centering
    \includegraphics[width=0.49\textwidth]{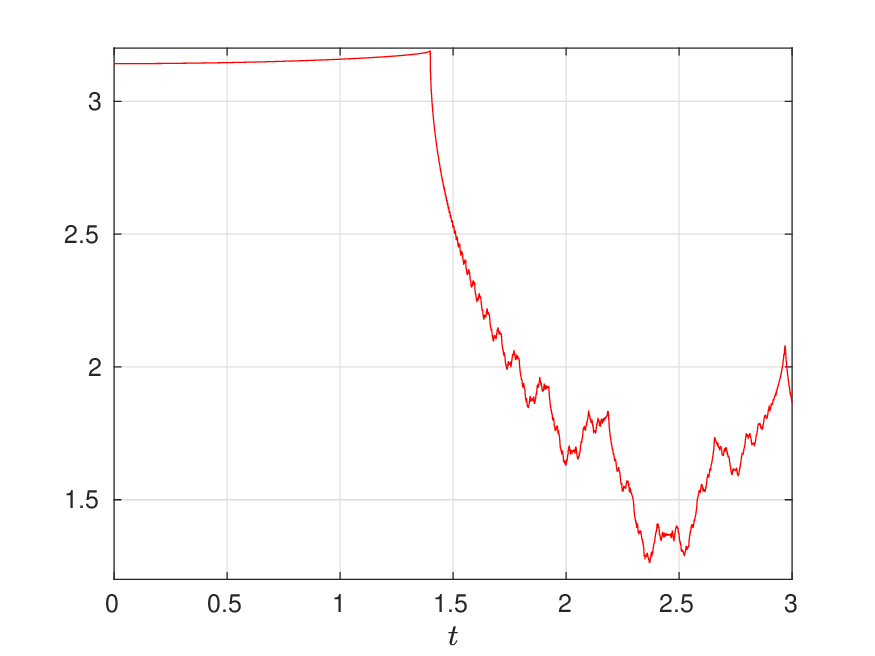}        
    \includegraphics[width=0.49\textwidth]{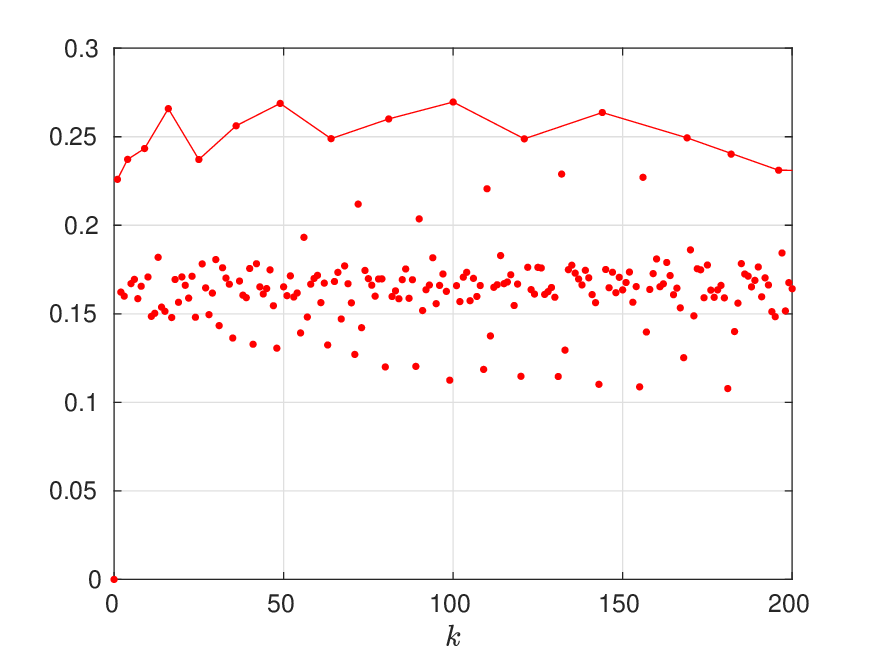}                    
    \caption{The modulus of the trajectory of the reconnection point $|\mathbf{X}(0,t)|$ (left) and its Fourier coefficients
    $k \widehat{|\mathbf{X}(0,t)|}(k)$ (right), calculated for the interval $(t_{rec}, t_{rec} + \pi / 2)$. 
    The line connects the coefficients corresponding to squares of integers.}\label{fig:rec_traj}
\end{figure}

Since the bridge is not present anymore ,we can follow the trajectory of the corner. It is possible to see that 
the frequencies corresponding to squares are dominating in the Fourier transform of the modulus of the trajectory (Figure~\ref{fig:rec_traj}). 
The transform here is performed in the interval $(t_{rec}, t_{rec} + \pi / 2)$ corresponding to a half of the quasi-period. Moreover,
the behavior of the reconnection point changes dramatically when the reconnection happens: from smooth to highly oscillating. It is reminiscent of
what happens in real fluids, where the reconnection of vortices may turn a laminar flow into a turbulent one.

\begin{figure}
    \centering
    \includegraphics[width=0.49\textwidth]{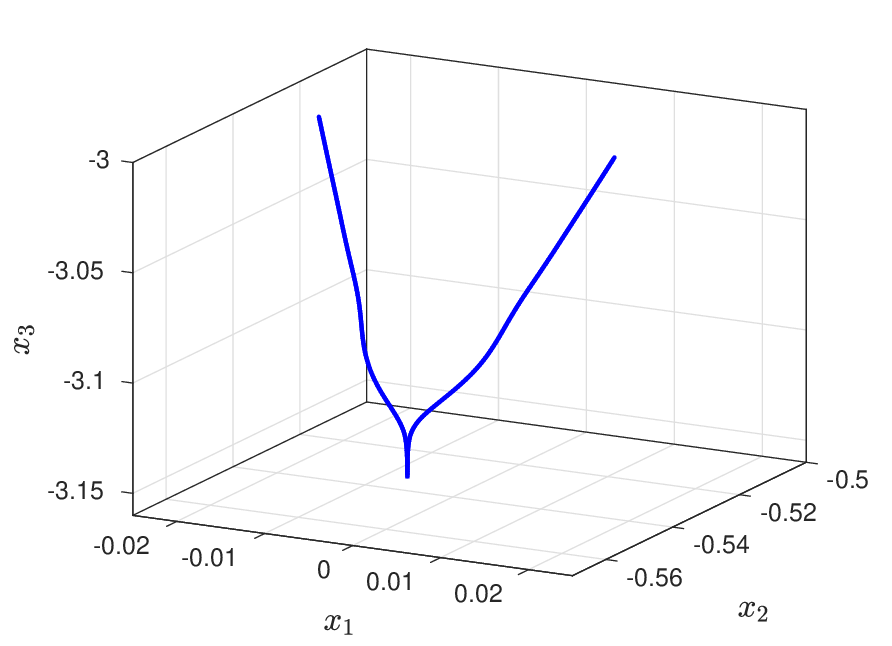}        
    \includegraphics[width=0.49\textwidth]{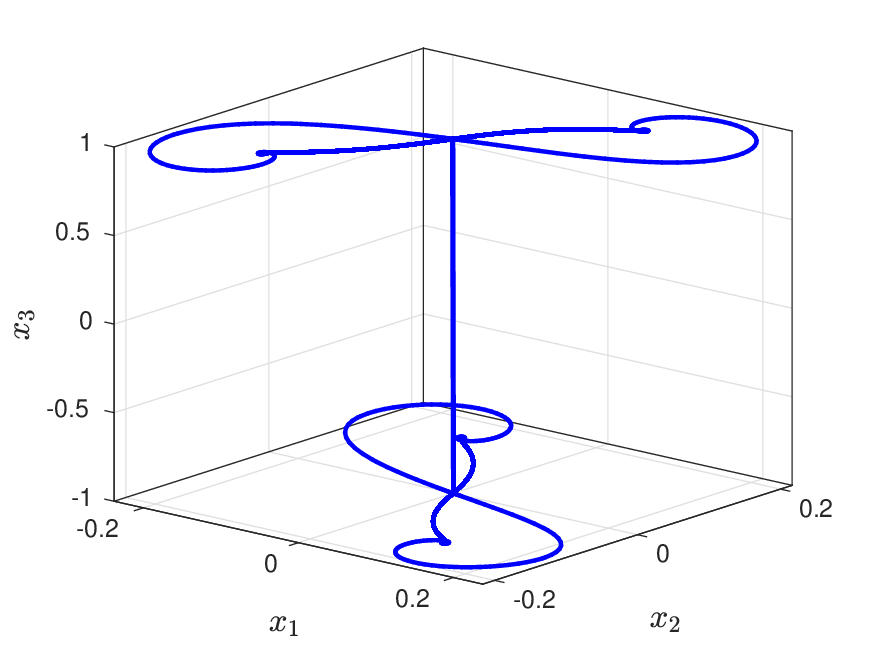}        
    \caption{Configuration of the vortices at the reconnection time $t_{rec}$, when $\min_s x_1(s) = 0$.} \label{fig:rec_moment}
\end{figure}

Despite the similarities, the shape of vortices at the reconnection time differs from the eye-shaped vortex~\eqref{eq:eye_vortex}; this is depicted in Figure~\ref{fig:rec_moment}. One can notice that there is a cusp instead of a corner. Moreover, the first horseshoe
has already appeared (left), as well as the first helical waves, which can be seen in the plot of the tangent vector (right). 
On the right plot, the vertical line is a jump between the values of the tangent vector around the cusp. The nodes that can be seen above and below are the first helical waves. 
After a short time, they turn into spirals as on the left-hand side of Figure~\ref{fig:rec_sol_tan}.  This supports the result in Section~\ref{sec:eye-shaped} that states that the behavior of the eye-shaped vortex does not depend on the angle, but it is dictated by the symmetry.

\subsubsection*{Analysis of the proposed reconnection procedure}

The results that we have for the selected values of $r_c$, $\varepsilon$, and the discretization considered are not so easily obtained when choosing arbitrary values. There are two reasons for it. 
Firstly, if the parameter $r_c$ is too big, then $x_1$ never reaches zero; therefore, we cannot use the reconnection criterion~\eqref{eq:min_u_crit}. 
Nevertheless, it is still possible to see the coherent structures, but with fewer details. This effect is reminiscent of what we can see in 
the solution of the Navier-Stokes equations for thick vortices. One possible way to perform the reconnection in this case is to consider another 
criterion, e.g., the fluid impulse. It was suggested in the numerical simulation in~\cite{iakunin2023} that the fluid impulse around the reconnection point is
monotone before the reconnection and oscillates afterwards. We can use this change of behavior, in order to define the reconnection time. Thus, for a 
threshold $\text{th}_{F}$, the reconnection time $t_{rec}$ is the smallest time when
\begin{equation}
    \frac{\partial}{\partial t}|\mathbf{F}_l(t_{rec})| \frac{\partial}{\partial t}|\mathbf{F}_l(t_{rec} - \tau)| \le -\text{th}_{F},  \label{eq:fi_crit}
\end{equation}
where $\tau$ is the time step. This allows us to perform the reconnection even for bigger values of $r_c$ that represent thicker vortices;
however, the smallest value of $x_1$ may not be zero. This may create an asymmetry that perturbs the solution and makes impossible to extract the frequencies corresponding to the squares.


The second reason why the reconnection of vortices is challenging for arbitrary values is related to numerical stability. 
Similarly as in the previous work~\cite{iakunin2023}, we use an embedded 5th Runge-Kutta method~\cite{butcher} in time with 8th order finite 
difference discretization in the filament parameter $s$. The necessary stability condition for the time-step $\tau$ and
spacial discretization step $h$ reads
\begin{equation}
    \tau \le \frac{h^2}{\sqrt{4 + \varepsilon  \frac{h^2}{r_c^2}}}.
\end{equation}
This means that we have to keep $h$ being of the order of $r_c$, in order to make the numerical method stable. Thus, it requires a lot of computational effort to model the reconnection for small values of $r_c$.


The reconnection procedure considered makes a lot of assumptions, e.g., the possibility of reparametrization of the vortex or the neglection of the
interaction. The reparametrization by arc-length has very little influence on the behavior,
because the vortex stretching at the reconnection time is not strong enough, whereas the interaction is more important, taking into account the second 
reconnection. We can notice especially that the interaction affects the quasi-period due to the change of the vortex length (formula~\eqref{eq:arclength}).

\begin{figure}
    \centering
    \includegraphics[width=0.49\textwidth]{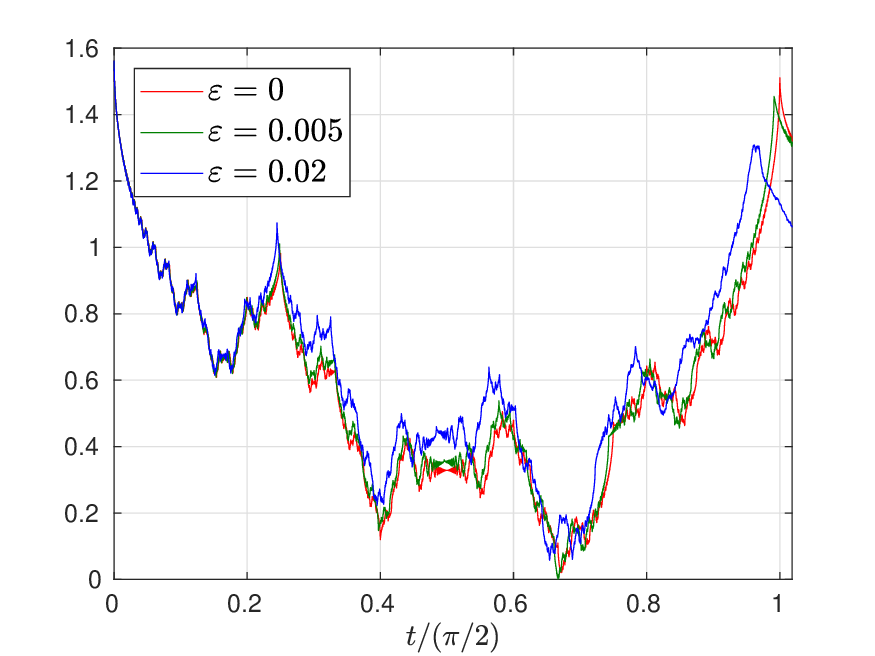}        
    \includegraphics[width=0.49\textwidth]{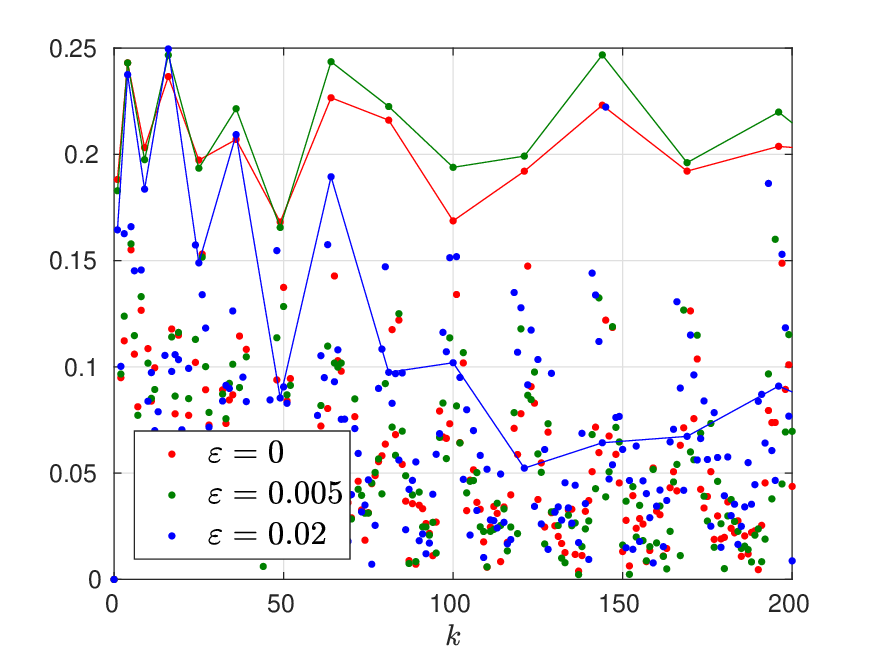}        
    \caption{The modulus of trajectory of the corner of the eye-shaped vortex $|\mathbf{X}(0,t)|$ (left) and its Fourier coefficients
    $k \widehat{|\mathbf{X}(0,t)|}(k)$ (right). 
    The line connects coefficients corresponding to squares of integers.} \label{fig:eps_eye_traj} 
\end{figure}

We can study how the interaction in~\eqref{eq:vfe_interaction} affects the frequencies in the toy model for the eye-shaped
vortex presented in section~\ref{sec:eye-shaped}, supposing now that this vortex obeys~\eqref{eq:vfe_interaction} for
some value of $\varepsilon > 0$. The results are depicted in Figure~\ref{fig:eps_eye_traj}. We can see that the quasi-period
is decreasing due to the vortex stretching. Furthermore, the stretching rate for the period is growing with time. It leads to
mixing the frequencies that correspond to squares with the remaining ones. For small values of $\varepsilon = 0.005$, the
frequencies that correspond to squares are still dominating; however, for $\varepsilon = 0.005$, they rapidly decrease. One possible 
direction for further research is finding a method to filter this stretching effect, using the expression~\eqref{eq:arclength} for the length of 
the tangent vector.

\section{Numerical simulation of the vortex reconnection in OpenFOAM} \label{sec:navier-stokes}

This section is devoted to the simulation of reconnection of vortices by the numerical solution of the incompressible Navier-Stokes equations (NSE)
\begin{align}
    &\mathbf{u}_t + \left(\mathbf{u} \cdot \nabla\right) \mathbf{u} = \nu \nabla^2 \mathbf{u} - \frac{1}{\rho} \nabla p, \text{ in } \Omega, \label{eq:ns_1}\\
    &\nabla \cdot \mathbf{u} = 0, \text{ in } \Omega,\label{eq:ns_2} 
\end{align}
in OpenFOAM. Here, $\mathbf{u}$ is the velocity of the fluid, $p$ is the pressure, $\rho$ is the density, $\nu = \mu / \rho$ is the kinetic viscocity, $\Omega \subset \mathbb{R}^3$
is a rectangular domain with periodic boundary conditions, and
$\nabla = \begin{pmatrix} \frac{\partial}{\partial x} & \frac{\partial}{\partial y} & \frac{\partial}{\partial z}\end{pmatrix}$.
We consider a domain $\Omega=(0,60)\times(0,144)\times(0,320)\ m^3$, filled with air with density $\rho = 1.2250\ kg / m^3$ 
and kinetic viscosity $\nu = 1.4607 \times 10^{-5} m^2 / s$. The initial condition is a pair of antiparallel Burnhamm-Hallock vortices:
\begin{equation}
    \mathbf{u}(\mathbf{x}, 0) = \frac{\Gamma}{2\pi} \left(\frac{\rho_1}{\rho_1^2 + r_c^2} (\mathbf{x} - \mathbf{x}_1) \wedge \mathbf{T} - 
                                           \frac{\rho_2}{\rho_2^2 + r_c^2} (\mathbf{x} - \mathbf{x}_2) \wedge \mathbf{T}\right),
\end{equation}
with core radius $r_c = 1.5\ m$, circulation $\Gamma = 3000\ m^2 / s$. $\rho_k$, with $k \in \{1,2\}$, is the closest distance from $\mathbf{x}$ to the $k$th vortex whose central line
is $\mathbf{x}_k + s\mathbf{T}$. This pair of vortices are perturbed by a random noise everywhere in the domain. 

\begin{figure}
    \includegraphics[width=\textwidth]{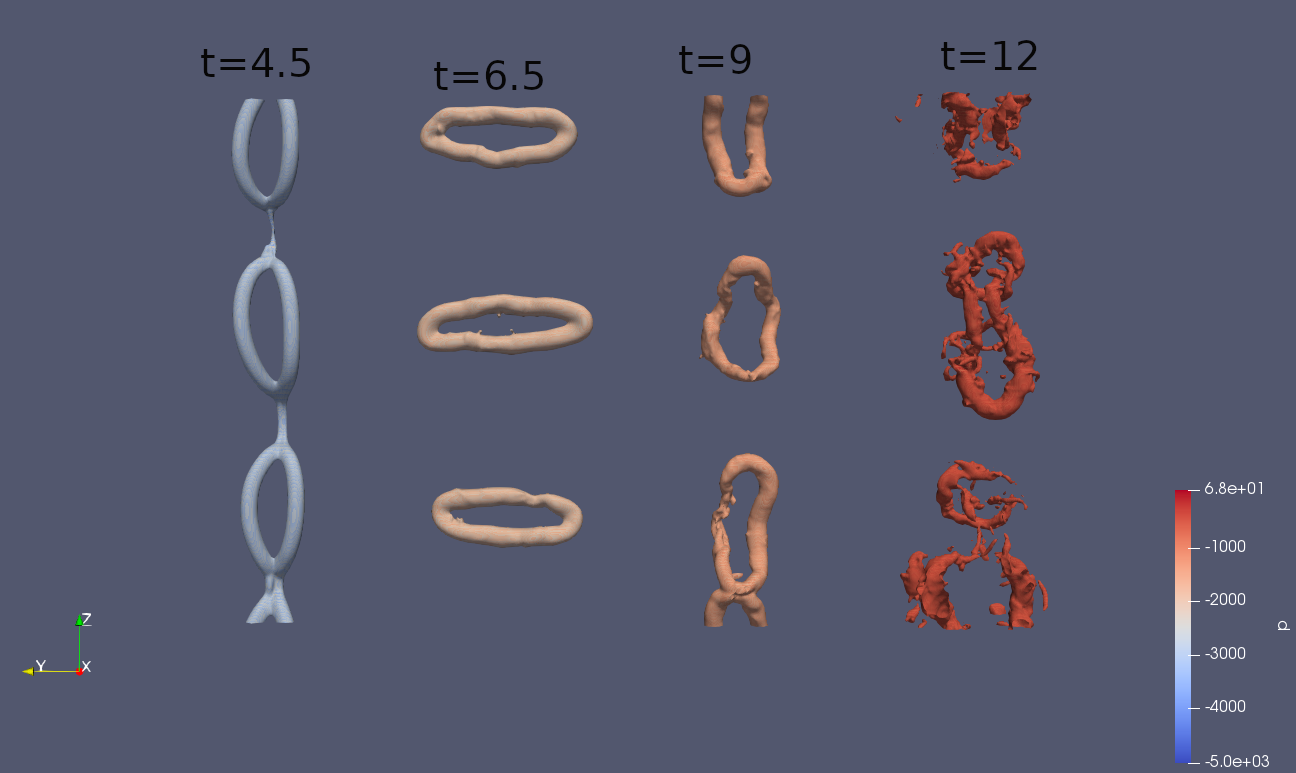}
    \caption{Level surface of the pressure field for the numerical solution of the NSE at different time moments.}\label{fig:ns_solution}
\end{figure}

We solve this problem numerically using the PIMPLE algorithm from OpenFOAM on a uniform mesh with size of $100\times 200 \times 400$ finite volumes. 
In Figure~\ref{fig:ns_solution}, the vortices are visualized at different times using pressure field. 
The reconnection happens at a time around $t = 4.5 s$. Then, at time $t = 6.5 s$, the configuration of vortices
looks also like an eye, but it is rotated. At time $t = 9s$, we have also eye-shaped vortices with similar orientation as at the reconnection time.
This process, however, does not repeat in time, because the vortices start to dissipate due to viscosity and the second reconnection. 

\begin{figure}
    \includegraphics[width=\textwidth]{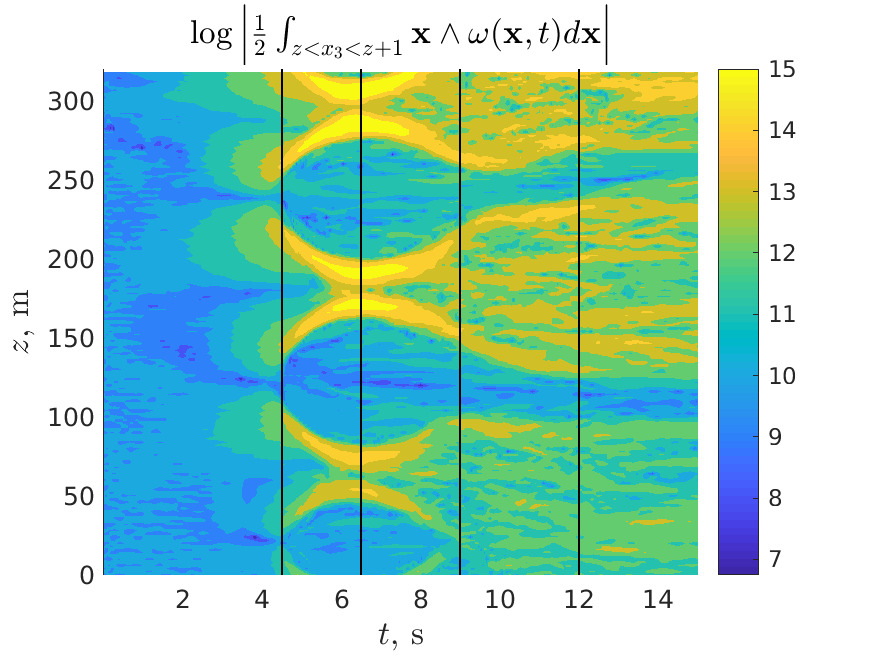}
    \caption{Modulus of the fluid impulse calculated for the numerical solution of the Navier-Stokes equations. 
    The lines correspond to the time moments shown in Figure~\ref{fig:ns_solution}}\label{fig:ns_fi}
\end{figure}

Even though it is quite challenging to extract the trajectory of the vortex central line from the numerical solution of NSE, the fluid impulse defined 
in~\eqref{eq:fluid_impulse} for infinitely thin vortices can be easily generalized for vortices of finite thickness. Consider a subdomain $\mathcal{L} \subset \Omega$;
then, the fluid impulse is
\begin{equation}
    \mathbf{F}_\mathcal{L}(t) = \frac{1}{2}\int_\mathcal{L} \mathbf{x}\wedge \boldsymbol{\omega}(\mathbf{x},t) d\mathbf{x},\label{eq:fi_ns}
\end{equation}
where $\boldsymbol{\omega} = \nabla \wedge \mathbf{u}$ is the vorticity. Similarly to the infinitely thin case, we are interested in the region $\mathcal{L}$ where
the reconnection happens. This region is not known from the initial conditions, but we can discretize the whole domain $\Omega$ in the $z$-direction with step $\delta z = 1\,m$
and define a function:
\begin{equation}
    F(z,t) = \left|\int_{0}^{L_x}\left(\int_{0}^{L_y}\left(\int_{z}^{z+\delta z}  \mathbf{x} \wedge \boldsymbol{\omega}(x,y,\zeta,t) d\zeta\right) dy \right) dx\right|,\label{eq:ns_fluid_impulse}
\end{equation}
where $L_x$ and $L_y$ are the size of $\Omega$ in the $x$- and $y$-directions respectively, with $\mathbf{x}=\begin{pmatrix} x & y & \zeta \end{pmatrix}$. 
The value of~\eqref{eq:ns_fluid_impulse} calculated for
the Navier-Stokes simulation is depicted in Figure~\ref{fig:ns_fi}. It is possible to see that, around the time $t = 4.5s$, the behavior of $F(z,t)$ changes from
homogeneous to oscillatory, with regions of high and low values. This effect is similar to what we have for the infinitely thin approximation. Furthermore,
we can clearly see periodicity in space that corresponds to Crow waves. The quasi-periodicity in time is also present. Note that the half of
the quasi-period corresponds to $t = 6.5s$, whereas its end is at $t = 9s$. Thus, the period is changing due to the viscosity and the vortex stretching,
similarly as in the case of infinitely thin vortices with the interaction term \eqref{eq:vfe_interaction}. The behavior for times $t > 9s$ is different, due to the dissipation caused by viscosity.

\begin{figure}
    \includegraphics[width=0.5\textwidth, height=0.4\textwidth]{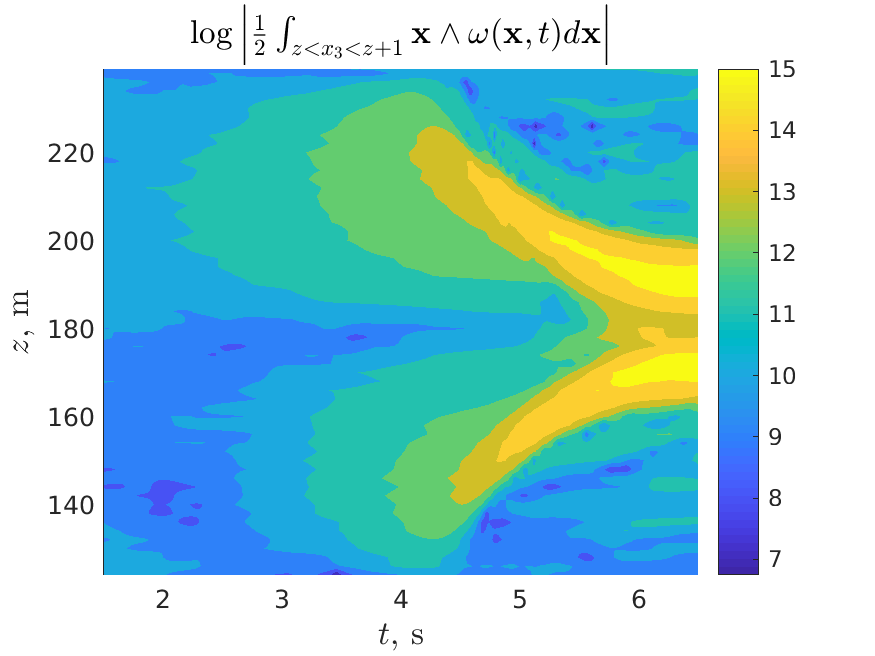}
    \includegraphics[width=0.5\textwidth, height=0.4\textwidth]{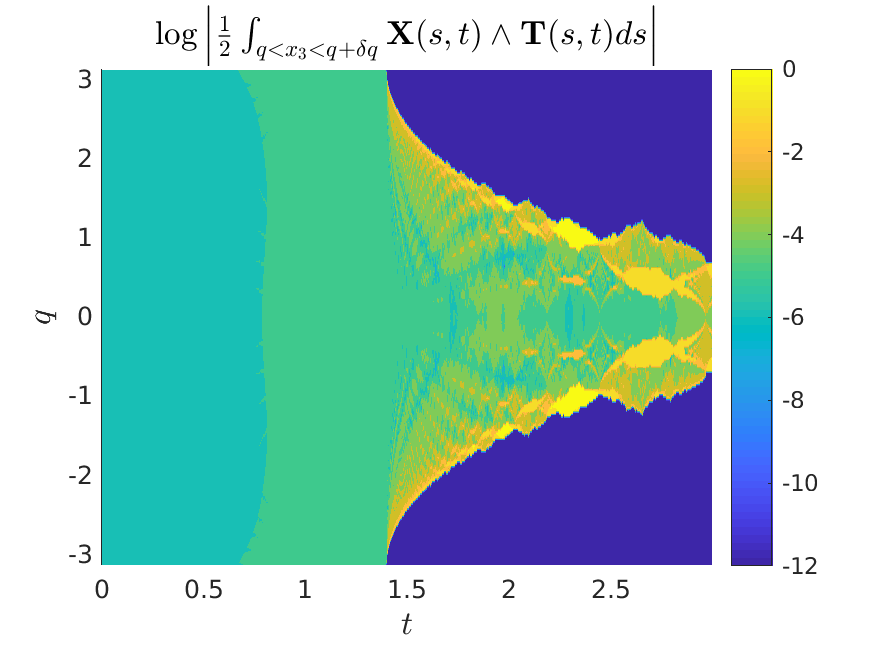}
    \caption{Comparison between the modulus of the fluid impulse of a finite thickness vortex (left) and an infinitely thin vortex (right)}\label{fig:comparison}
\end{figure}

We can apply \eqref{eq:ns_fluid_impulse} to the solution of~\eqref{eq:vfe_interaction} considered in the previous paragraph. In this case, we introduce
a parameter $q$, a small step $\delta q$ and calculate the following integral:
\begin{equation}
    F(q,t) = \left|\int_{q < x_3(s,t) < q + \delta q} \mathbf{X}(s,t) \wedge \mathbf{T}(s,t)ds\right|. \label{eq:ns_fluid_impulse_filament}
\end{equation}
The comparison between~\eqref{eq:ns_fluid_impulse} and~\eqref{eq:ns_fluid_impulse_filament} is depicted in Figure~\ref{fig:comparison}. One can see that
in the infinitely thin case much more details are present. Note also that the profile of non-zero values after the reconnection time ($t_{rec} = 1.397$) on the right-hand side of 
Figure~\ref{fig:comparison} is very reminiscent of the modulus of the corner trajectory (Figure~\ref{fig:rec_traj}) or RNDF~\eqref{eq:rndf}. It indeed represents the third
component of the corner trajectory and can be used even for a vortex with finite thickness. In particular, on the left-hand side of Figure~\ref{fig:comparison}, after
the reconnection time ($t = 4.5s$), there are regions of fluid impulse concentration whose profile is approximately $\sqrt{t}$. This is reminiscent of the vortex separation 
rate observed in quantum fluids~\cite{fonda2019}. This rate was also obtained in the previous work~\cite{iakunin2023} for infinitely thin vortices, and derived analytically for the 
vortex made by two intersected half-lines~\cite{vega2003}. There are fewer oscillations than in the case of an infinitely thin vortex, but, if we perform a finer simulation of NSE~\eqref{eq:ns_1},~\eqref{eq:ns_2}, we can expect a more detailed result.

\begin{figure}
    \centering 
    \includegraphics[width=\textwidth]{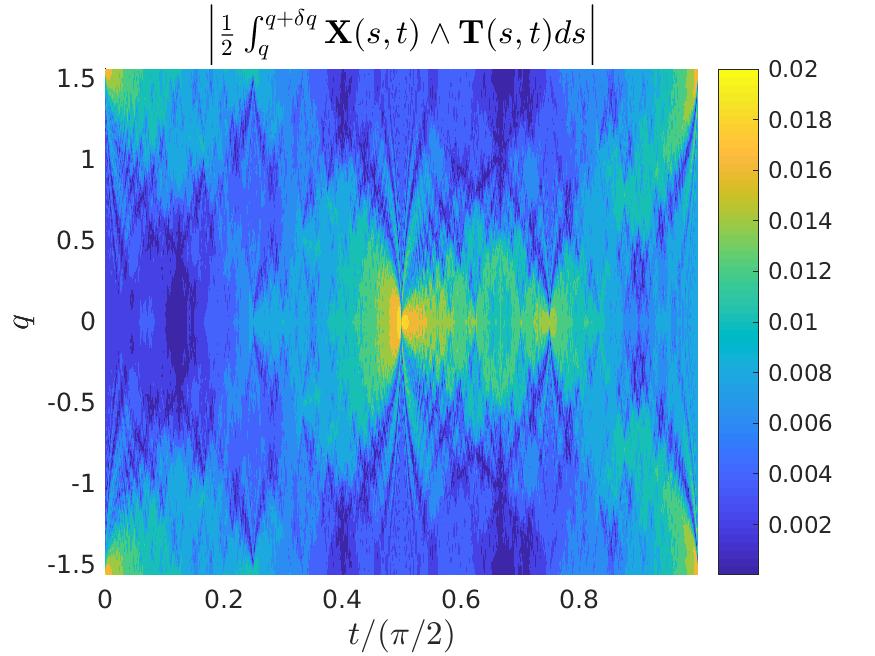}
    \caption{Modulus of the fluid impulse calculated for the eye-shaped vortex}\label{fig:vfe_eye_fi_surf}
\end{figure}

Another effect that can be noticed is the emergence of rhombus-like structures in the fluid impulse plot (right-hand side of Figure~\ref{fig:comparison}) after the reconnection. These
structures can be obtained in a clearer way if we perform the integration of \eqref{eq:ns_fluid_impulse_filament} along the parameter of the curve:
\begin{equation}
    F(q,t) = \left|\int_{q}^{q + \delta q} \mathbf{X}(s,t) \wedge \mathbf{T}(s,t)ds\right|. \label{eq:ns_fluid_impulse_filament_s}
\end{equation}
The result for the eye-shaped vortex~\eqref{eq:eye_vortex} is depicted in Figure~\ref{fig:vfe_eye_fi_surf}. It resembles the Talbot effect in quantum optics, i.e., the interference of particles passing through holes on a screen. The corners of the eye-shaped vortex at time $t=0$ are zones with high intensity of $F(q,t)$, and are similar to the holes in the Talbot effect. At time $t = \pi / 4$, which corresponds to the half of a quasi-period, the fluid impulse in concentrated in the center, which
is also similar to the Talbot effect. At the end, at time $t = \pi / 2$, the zones of high intensity of $F(q,t)$ are the same as they were initially, which is also present
in the Talbot effect. The relation between the Talbot effect and NLSE, which can be obtained from VFE by the Hasimoto transform, is studied by many authors, e.g., in~\cite{erdogan2013}.
Its presence in the fluid impulse, however, suggests the possibility of observing this effect, or some approximation of it, in the solution of NSE. 

\section{Conclusions} \label{sec:conclusions}

The phenomenon of the reconnection of a pair of antiparallel vortices is studied in both the finite thickness case and the
infinitely thin approximation case. One of the main questions one can ponder is if the vortices form a corner, and
what is the definition of the reconnection time and the reconnection point. It is known that the behavior of a corner
singularity of an infinitely thin closed vortex is very complex and possibly multi-fractal. Therefore, the capacity
to discover such behavior in the vortex reconnection can be beneficial for understanding the turbulence. In this article, we propose
a model for the reconnection of infinitely thin vortices that allows to observe a behavior reminiscent of the one of a vortex in the shape of a regular polygon.
We also show that the fluid impulse can be used for the definition of the reconnection time and the reconnection point, even for the
finite thickness case.

We have performed the study of an isolated eye-shaped vortex as an approximation of the configuration of vortices after the reconnection. We have observed the existence of quasi-periodicity independent of the angle of the corner, in the sense that an eye-shaped vortex with the same orientation (but eventually with different values of the angle of the corner) is recovered after a quasi-period in time of value $\pi/2$. This unexpected result allows to generalize the reconnection process, since the configuration of vortices after the reconnection is always close to an eye with different corner angles. Moreover, we have observed that the conservation of the total of the fluid impulse might be used in order to discover more preserved quantities.

The results obtained for eye-shaped vortices enable us develop a heuristic approach for the reconnection of infinitely thin vortices, which allows to get rid of the bridge. The behavior of vortices after the reconnection appears to be very reminiscent of the behavior of vortices in the shape of a regular polygon or in the shape of an eye. In particular, we are able to observe in the Fourier transform of the trajectory of the reconnection point the domination of frequencies that correspond to squares of integers. This can suggest that the behavior of a vortex after reconnection may be related to Riemann's non-differentiable function, and hence exhibit multi-fractal properties.

However, the application of the reconnection approach developed is delicate, and depends on the parameters of the vortex reconnection model, such as $\varepsilon$ and $r_c$. This is due to the fact that, even though the shape of a vortex at the reconnection time is close to an eye, it does not exactly resemble it. In particular, we can observe cusps instead of corners, and a small horseshoe structure, whose size depends on $r_c$.  In order to get rid of this effect, we have to set $r_c = 0$, which is not possible, due to numerical stability restrictions.

We have also studied how the interaction between vortices influences the multi-fractal properties of the trajectory of the corner of an eye-shaped vortex. In this regard, we have observed that the speed of the process increases, due to the stretching of the vortex. Thus, there is no more a quasi-period, and the frequencies of the Fourier transform are mixed. One way to avoid it is to make corrections, including the vortex stretching rate, which can be obtained from the expression for the tangent vector \eqref{eq:arclength}.

Some features of the behavior of the infinitely thin vortices can be also observed in the solution of the Navier-Stokes equations, whose initial condition is a pair of antiparallel vortices of finite thickness. We have considered the modulus of the fluid impulse, and noticed that its behavior changes from monotone before the reconnection to oscillatory after the reconnection, for both the finite thickness case, and the infinitely thin case. Moreover, the first quasi-period can be observed, whereas the successive quasi-periods are no longer distinguishable, due to viscosity. The amount of details that can be appreciated in the case of a finite thickness vortex is much smaller that in the infinitely thin case, but we can expect more resemblances for finer simulations.





    \appendix

\section{Derivation of the conserved quantity related to the symmetry}\label{sec:first_cons_quant}

Due to the symmetry of the rhombus, the product $\mathbf M_1\mathbf M_0$ induces a rotation of $\pi$ radians, so its trace is $1 + 2\cos(\pi) = -1$. On the other hand, it is given by:
\begin{align*}
    \operatorname{trace}{M_1 M_0} &= \cos{\rho_0}\cos{\rho_1} \cr
    &+ \left(\cos{\rho_0} \cos^2{\theta_0} - \cos^2{\theta_0} + 1\right) \left(\cos{\rho_1} \cos^2{\theta_1} - \cos^2{\theta_1} + 1\right) \cr
    &+ \left(\cos^2{\theta_0} + \cos{\rho_0}\sin^2{\theta_0}\right) \left(\cos^2{\theta_1} + \cos{\rho_1}\sin^2{\theta_1}\right) \cr
    &+ 2\sin{\theta_0} \cos{\theta_0} \sin{\theta_1} \cos{\theta_1} \left(\cos{\rho_0} - 1\right) \left(\cos{\rho_1} - 1\right) \cr
    &- 2 \cos{\theta_0} \cos{\theta_1} \sin{\rho_0}\sin{\rho_1} - 2\sin{\rho_0}\sin{\rho_1}\sin{\theta_0}\sin{\theta_1}.
\end{align*}
The term independent of $\rho_0$ and $\rho_1$ reads:
\begin{align*}
    & \left(1-\cos^2{\theta_0}\right) \left(1-\cos^2{\theta_1}\right) + \cos^2{\theta_0} \cos^2{\theta_1} + 2\sin{\theta_0} \cos{\theta_0} \sin{\theta_1} \cos{\theta_1}
    \cr
    & \qquad = \cos^2(\theta_0 - \theta_1).
\end{align*}
Here, we have used that $\cos{\theta_0} \cos{\theta_1} + \sin{\theta_0} \sin{\theta_1} = \cos(\theta_0 - \theta_1)$ and the Pythagorean identity. 
The terms proportional to $\cos{\rho_0}$ and $\cos{\rho_1}$ are respectively
\begin{equation*}
    \cos^2{\theta_0}\left(1-\cos^2{\theta_1}\right) + \sin^2{\theta_0}\cos^2{\theta_1} - 2\sin{\theta_0} \cos{\theta_0} \sin{\theta_1} \cos{\theta_1} = \sin^2(\theta_0 - \theta_1),
\end{equation*}
and
\begin{equation*}
    \left(1-\cos^2{\theta_0}\right)\cos^2{\theta_1} + \sin^2{\theta_1}\cos^2{\theta_0} - 2\sin{\theta_0} \cos{\theta_0} \sin{\theta_1} \cos{\theta_1} = \sin^2(\theta_0 - \theta_1).
\end{equation*}
Here, we have used that $\sin{\theta_0} \cos{\theta_1} - \cos{\theta_0} \sin{\theta_1} = \sin(\theta_0 - \theta_1)$. Finally, the term proportional to $\cos{\rho_0} \cos{\rho_1}$ is
\begin{equation*}
    1 + \cos^2{\theta_0}\cos^2{\theta_1} + \sin^2{\theta_0} \sin^2{\theta_1} + 2\sin{\theta_0} \cos{\theta_0} \sin{\theta_1} \cos{\theta_1} = 1 + \cos^2(\theta_0 - \theta_1),
\end{equation*}
and the term proportional to $\sin{\rho_0} \sin{\rho_1}$ is
\begin{equation*}
    - 2 \cos{\theta_0} \cos{\theta_1} - 2\sin{\rho_0}\sin{\rho_1} = -2\cos(\theta_0 - \theta_1).
\end{equation*}
Thus, we can simplify the expression for the trace:
\begin{align*}
    \operatorname{trace}{M_1 M_0} &= \cos^2(\theta_0 - \theta_1) \cr
    &\quad+ \left(\cos{\rho_0} + \cos{\rho_1}\right) \sin^2(\theta_0 - \theta_1) \cr
    &\quad+ \cos{\rho_0} \cos{\rho_1} \left(1 + \cos^2(\theta_0 - \theta_1)\right) \cr
    &\quad-2 \sin{\rho_0} \sin{\rho_1}\cos(\theta_0 - \theta_1).
\end{align*}
Moreover, we can get rid of $\sin^2(\theta_0 - \theta_1)$ by means of the Pythagorean identity:
\begin{align*}
    \operatorname{trace}{M_1 M_0} &= \cos{\rho_0} + \cos{\rho_1} + \cos{\rho_0} \cos{\rho_1} \cr
    &\quad- 2 \sin{\rho_0} \sin{\rho_1}\cos(\theta_0 - \theta_1) \cr
    &\quad+  \left(1 - \cos{\rho_0}\right)\left(1 - \cos{\rho_1}\right)  \cos^2(\theta_0 - \theta_1).
\end{align*}
Note that, for any $\alpha$, we have that $\cos{\alpha} = 2\cos^2(\alpha/2) - 1 = 1 - 2\sin^2(\alpha/2)$ and $\sin{\alpha} = 2 \cos(\alpha/2)\sin(\alpha/2)$. 
Applying these formulas to $\rho_0$ and $\rho_1$, we obtain:
\begin{align*}
    \operatorname{trace}{M_1 M_0} &= -1 \cr
    &\quad+ 4 \cos^2(\rho_0 / 2) \cos^2(\rho_1 / 2) \cr
    &\quad- 8 \sin(\rho_0 / 2) \cos(\rho_0 / 2) \sin(\rho_1 / 2) \cos(\rho_1 / 2) \cos(\theta_0 - \theta_1) \cr
    &\quad+ 4 \sin^2(\rho_0 / 2) \sin^2(\rho_1 / 2) \cos^2(\theta_0 - \theta_1),
\end{align*}
which simplified becomes
\begin{align*}
   \operatorname{trace}{M_1 M_0} & = -1
   \cr
   & \quad + 4 \sin^2(\rho_0 / 2) \sin^2(\rho_1 / 2) (\cos(\theta_0 - \theta_1)- \cot(\rho_0 / 2) \cot(\rho_1 / 2))^2.    
\end{align*}
Due to the symmetry, the trace is equal to $-1$, which yields \eqref{eq:trace_conservation}.

\section{Derivation of the conserved quantity related to the fluid impulse}\label{sec:second_cons_quant}
Let us consider the length of a vector $\mathbf{f} = \mathbf{T}_0 \wedge \mathbf{T}_1 + \mathbf{T}_2 \wedge \mathbf{T}_3$, such that the total
fluid impulse if given by $\mathbf{F} = \pi^2 / 8\ \mathbf{f}$. 
For any four vectors $\mathbf{a}$, $\mathbf{b}$, $\mathbf{c}$, $\mathbf{d}$, it is known that
\begin{equation*}
    (\mathbf{a} \wedge \mathbf{b})\cdot(\mathbf{c} \wedge \mathbf{d}) = (\mathbf{a}\cdot\mathbf{c})(\mathbf{b}\cdot\mathbf{d}) - (\mathbf{a}\cdot\mathbf{d})(\mathbf{b}\cdot\mathbf{c}).
\end{equation*}
Thus, the square of the length of $\mathbf{f}$ is 
\begin{align}
    |\mathbf{f}|^2 &= (\mathbf{T}_0 \wedge \mathbf{T}_1)^2 + (\mathbf{T}_2 \wedge \mathbf{T}_3)^2 + 2 (\mathbf{T}_0 \wedge \mathbf{T}_1)\cdot (\mathbf{T}_2 \wedge \mathbf{T}_3) \cr
    &= |\mathbf{T}_0|^2 |\mathbf{T}_1|^2 - (\mathbf{T}_0 \cdot \mathbf{T}_1)^2 + |\mathbf{T}_2|^2 |\mathbf{T}_3|^2 - (\mathbf{T}_2 \cdot \mathbf{T}_3)^2 \cr
    &\quad+ 2 (\mathbf{T}_0 \cdot \mathbf{T}_2) (\mathbf{T}_1 \cdot \mathbf{T}_3) - 2 (\mathbf{T}_0 \cdot \mathbf{T}_3) (\mathbf{T}_1 \cdot \mathbf{T}_2) \cr
    &= 2-(\mathbf{T}_0\cdot\mathbf{T}_1)^2-(\mathbf{T}_2\cdot\mathbf{T}_3)^2
    \cr
    & \quad +2 (\mathbf{T}_0\cdot\mathbf{T}_2) (\mathbf{T}_1\cdot\mathbf{T}_3)-2 (\mathbf{T}_0\cdot\mathbf{T}_3) (\mathbf{T}_1\cdot\mathbf{T}_2).\label{eq:fi_simplified}
\end{align}
Here, we have used that the length of all tangent vectors is $1$. Using the symmetry, we can conclude 
that $\mathbf{T}_0 \cdot \mathbf{T}_1 = \mathbf{T}_2 \cdot \mathbf{T}_3 = \cos{\rho_0}$, as well as 
$\mathbf{T}_0 \cdot \mathbf{T}_3 = \mathbf{T}_1 \cdot \mathbf{T}_2 = \cos{\rho_1}$. In order to prove 
$\mathbf{T}_0 \cdot \mathbf{T}_3 = \cos{\rho_1}$ directly, we have to use the property~\eqref{eq:trace_conservation} 
proved in Appendix~\ref{sec:first_cons_quant}. The last unknown scalar products are $\mathbf{T}_0 \cdot \mathbf{T}_2$ and $\mathbf{T}_1 \cdot \mathbf{T}_3$, which are the same due to the symmetry. From the definition of $\mathbf{T}_0$ and $\mathbf{T}_2$, we have:
\begin{align*}
    \mathbf{T}_0 \cdot \mathbf{T}_2 &= \cos{\rho_0} \cos{\rho_1} - \cos{\theta_0} \cos{\theta_1}\sin{\rho_0}\sin{\rho_1} - \sin{\rho_0}\sin{\rho_1} \sin{\theta_0}\sin{\theta_1}\cr
    &=\cos{\rho_0} \cos{\rho_1} - \cos(\theta_0 - \theta_0)\sin{\rho_0}\sin{\rho_1}\cr
    &=\left(2\cos^2(\rho_0 / 2) - 1\right)\left(2\cos^2(\rho_1 / 2)- 1\right) \cr
    & \quad - 4 \cot(\rho_0 / 2)\cot(\rho_1 / 2) \sin(\rho_0 / 2)\cos(\rho_0 / 2)\sin(\rho_1 / 2)\cos(\rho_1 / 2)\cr
    &=1 - 2\cos^2(\rho_0 / 2) - 2\cos^2(\rho_1 / 2) \cr
    &=-1 - \cos{\rho_0} - \cos{\rho_1}.
\end{align*}
Here, we have used the property \eqref{eq:trace_conservation} in order to get rid of $\cos(\theta_0 - \theta_0)$. Thus, the expression~\eqref{eq:fi_simplified} reads
\begin{align*}
    |\mathbf{f}|^2 &= 2 - 2\cos^2{\rho_0} + 2(1 + \cos{\rho_0} + \cos{\rho_1})^2 - 2 \cos^2{\rho_1}\cr
    &=4 + 4\cos{\rho_0} + 4\cos{\rho_1} + 4\cos{\rho_0}\cos{\rho_1}\cr
    &=4 (1 + \cos{\rho_0})(1 + \cos{\rho_1}),
\end{align*}
which is equivalent to~\eqref{eq:fi_conservation}. The expression~\eqref{eq:fi_expression} can be obtained without applying~\eqref{eq:trace_conservation} to~\eqref{eq:fi_simplified}.
In this case $\mathbf{T}_1 \cdot \mathbf{T}_2 = \cos{\rho_1}$ whereas 
\begin{align*}
    \mathbf{T}_0\cdot\mathbf{T}_3 & = -1+\cos^2{\rho_0}(\cos{\rho_1}+1)-2 \cos(\theta_0-\theta_1)\cos{\rho_0} \sin{\rho_0} \sin{\rho_1}
    \cr
    & \quad -
    \cos^2(\theta_0-\theta_1) \sin^2{\rho_0} (\cos{\rho_1} - 1).
\end{align*}    

    \section*{Acknowledgments}

Sergei Iakunin was partially supported by the BCAM Severo Ochoa accreditation CEX2021-001142-S. Francisco de la Hoz was partially supported by the research group grant IT1615-22 funded by the Basque Government. Francisco de la Hoz and Sergei Iakunin were partially supported by the grant PID2021-126813NB-I00 funded by MICIU/AEI/10.13039/501100011033 and by ``ERDF A way of making Europe''.	

The authors are grateful to Prof. Luis Vega, for useful discussions on this paper.

	\bibliographystyle{unsrt}

	\bibliography{bibliography}

\end{document}